\documentclass[11pt]{article}

\usepackage{amsmath}
\usepackage{amssymb}
\usepackage{enumerate}
\usepackage{graphicx}
\usepackage{pifont}
\usepackage{epsfig}
\usepackage{theorem}
\usepackage{xcolor}
\usepackage{algorithmic, algorithm}

\newcommand{\menge}[2]{\big\{{#1}~\big |~{#2}\big\}} 
\newcommand{\emp}{\ensuremath{\varnothing}}
\newcommand{\Frac}[2]{\displaystyle{\frac{#1}{#2}}} 
\newcommand{\scal}[2]{{#2}^\top {#1}} 
\newcommand{\HH}{\ensuremath{{\mathcal H}}}

\newcommand{\JJ}{\ensuremath{S}}
\newcommand{\II}{\ensuremath{\mathbb I}}

\newcommand{\RR}{\ensuremath{\mathbb R}}

\newcommand{\RPP}{\ensuremath{\,\left]0,+\infty\right[}}
\newcommand{\RX}{\ensuremath{\,\left]-\infty,+\infty\right]}}

\newcommand{\NN}{\ensuremath{\mathbb N}}

\newcommand{\dom}{\ensuremath{\mathrm{dom}\,}}
\newcommand{\Argmin}{\ensuremath{\mathrm{Argmin}\,}}

\newcommand{\prox}{\ensuremath{\mathrm{prox}}}

\newcommand{\rint}{\ensuremath{\mathrm{rint}\,}}
\newcommand{\ran}{\ensuremath{\mathrm{ran}\,}}

\newcommand{\Idf}{\ensuremath{\mathrm{Id}}}

\newcommand{\ic}{\ensuremath{\mathrm{\iota_C}}}

\newcommand{\pinf}{\ensuremath{+\infty}}

\newcommand{\sign}{\ensuremath{\mathrm{sign}}}
\newcommand{\diag}{\ensuremath{\mathrm{Diag}}}

\newcommand{\vv}{\ensuremath{V}}
\newcommand{\hh}{\ensuremath{H}}
\newcommand{\tv}{\ensuremath{\mathrm{tv}}}

\newcommand{\trace}{\ensuremath{\mathrm{tr}}}

\newtheorem{theorem}{Theorem}[section]
\newtheorem{proposition}[theorem]{Proposition}
\theoremstyle{plain}{\theorembodyfont{\rmfamily}%
\newtheorem{example}[theorem]{Example}}
\theoremstyle{plain}{\theorembodyfont{\rmfamily}%

\newtheorem{assumption}[theorem]{Assumption}

\newtheorem{remark}[theorem]{Remark}}

\title{Parallel Proximal Algorithm for Image Restoration Using Hybrid Regularization -- Extended Version \thanks{Part of this work appeared in the conference proceedings of EUSIPCO 2009 \cite{Pustelnik_N_2009_eusipco_Hybrid_rfdritpopn}. This work was supported by the Agence Nationale de la Recherche under grant ANR-09-EMER-004-03.}}

\author{Nelly Pustelnik, Caroline Chaux, and Jean-Christophe Pesquet
\thanks{N. Pustelnik (Corresponding Author), C. Chaux and J.-C. Pesquet are with the Universit{\'e} Paris-Est, Laboratoire d'Informatique Gaspard Monge, CNRS-UMR 8049, 77454 Marne-la-Vall{\'e}e
Cedex 2, France. Phone: +33 1 60 95 77 39, E-mail: \{nelly.pustelnik,caroline.chaux,jean-christophe.pesquet\}@univ-paris-est.fr}.
}

\begin{document}
\maketitle

\begin{abstract}
Regularization approaches have demonstrated their effectiveness
for solving ill-posed problems. However, in the context of variational restoration methods, a challenging question remains, namely how to find a good regularizer. While total variation introduces staircase effects, wavelet domain regularization brings other artefacts, e.g. ringing. However, a trade-off can be made by introducing a hybrid regularization including several terms non necessarily acting in the same domain (e.g. spatial and wavelet transform domains).
While this approach was shown to provide good results for
  solving deconvolution problems in the presence of additive Gaussian noise,
  an important issue is to efficiently deal with this
  hybrid regularization for more general noise models. To solve this problem, we adopt a convex optimization framework where the criterion to be minimized is split in the sum of more than two terms. For spatial domain regularization, isotropic or anisotropic total variation definitions using various gradient filters are considered.
An accelerated version of the Parallel Proximal Algorithm is proposed to perform the minimization. Some difficulties in the computation of the proximity operators involved in this algorithm are also addressed in this paper.
Numerical experiments performed in the context of Poisson data recovery, show the good behaviour of the algorithm as well as promising results concerning the use of hybrid regularization techniques.
\end{abstract}

\newpage
\section{Introduction}\label{sec:intro} 
During the last decades, convex optimization methods have been shown to be very effective for solving inverse problems.
On the one hand, algorithms such as Projection Onto Convex Sets (POCS) \cite{Bregman_LM_1965_sm_POCS_tmospfacpocs,Gurin_LG_1967_zvmmf_Projection_mffacpocs,Youla_DC_1982_tmi_POCS_irbtmopocs, Combettes_P_1993_pieee_fou_ste, Combettes_P_1996_book_con_fpi} 
have become popular for finding a solution in the intersection of convex sets. POCS was used in data recovery problems \cite{Trussell_H_1984_tassp_feasible_ssp} in order 
to incorporate prior information on the target image (e.g. smoothness  constraints). Some variants of POCS such as ART (Algebraic Reconstruction Technique) \cite{Gordon_R_1970_jtb_ART_artftdemaxrp} or PPM (Parallel Projection Method)  
\cite{DePierro_AR_1986_nm_PPM_crfaanca,Combettes_PL_1994_tsp_Inconsistent_sfplssiaps} 
were also proposed to achieve iteration parallelization. Additional variants of POCS can be found in \cite{Combettes_P_1997_tip_cn_sti}. Other parallel approaches such as block-iterative surrogate constraint splitting methods were considered to solve a quadratic minimization problem under convex constraints \cite{Combettes_PL_2003_tsp_Block_abiscsmfqsr} which may include a total variation constraint 
(see also \cite{Aujol_JF_2009_jmiv_firstorder_atvbir}). 
However the method in
\cite{Combettes_PL_2003_tsp_Block_abiscsmfqsr} based on subgradient projections is not applicable to non-differentiable objective functions.

On the other hand, some denoising approaches were based on wavelet transforms \cite{Mallat_S_1999_ap_wavelet_awtosp}, and more generally on frame representations 
\cite{Candes_EJ_2002_as_Recovering_eiipipoocf,LePennec_E_2005_tip_spa_girb,Selesnick_I_2005_dsp_dual_tdtcwt,Chaux_C_2006_tip_ima_adtmbwt}.
 In \cite{Daubechies_I_2004_cpamath_iterative_talipsc,Figueiredo_M_2003_tosp_EM_afwbir,Bect_J_2004_eccv_unified_vfir,Harmany_Z_2011_spiral_tap}, algorithms  which belong to the class of forward-backward algorithms were proposed in order to restore images 
degraded by a linear operator and a noise perturbation. Forward-backward iterations allow us to minimize 
a sum of two functions assumed to be in the class $\Gamma_0(\HH)$ of lower semicontinuous convex functions defined on a Hilbert space $\HH$, and taking their values in $]-\infty,+\infty]$, which are proper (i.e. not identically equal to $+\infty$). In addition,  one of these functions must be Lipschitz differentiable on $\HH$.
In \cite{Combettes_PL_2005_mms_Signal_rbpfbs}, this algorithm was investigated
by making use of proximity operator tools \cite{Combettes_P_2010_inbook_proximal_smsp} firstly proposed by Moreau in \cite{Moreau_J_1965_bsmf_Proximite_eddueh}. In \cite{Chaux_C_2007_ip_variational_ffbip}, applications to frame representations were developed and a list of closed form expressions of several proximity operators was provided. Typically,
forward-backward methods are appropriate when dealing with a smooth data fidelity
term e.g. a quadratic function and a non-smooth penalty term such as an 
$\ell_1$-norm promoting sparsity in the considered frame \cite{Bioucas_J_2007_toip_New_ttsistafir,Beck_A_2009_j-siam-is_fast_istalip}. The computation of the proximity operator associated with the $\ell_1$-norm indeed reduces to a componentwise soft-thresholding \cite{Donoho_D_1995_tit_den_st,Combettes_PL_2007_jopt_Proximal_tafmoob}. Another optimization method known as the Douglas-Rachford algorithm \cite{Lions_PL_1979_jna_Splitting_aftsotno,Douglas_J_1956_tams_numerical_otnsothcpi2a3sv,Eckstein_J_2003_mp_douglas_otdrsmatppaftmmo,Combettes_PL_2007_istsp_Douglas_rsatncvsr} was then proposed for signal/image recovery problems \cite{Combettes_PL_2007_istsp_Douglas_rsatncvsr} to relax the Lipschitz differentiablity condition required in forward-backward iterations. 
In turn, the latter algorithm requires the knowledge of the proximity operators of both functions. This algorithm was then extended to the minimization of 
a sum of a finite number of convex functions  \cite{Combettes_PL_2008_journal_proximal_apdmfscvip}, the proximity operator of each function still being assumed to be known. One of the main advantages of this algorithm called Parallel ProXimal Algorithm (PPXA) is its parallel
structure which makes it easily implementable on multicore architectures. PPXA is well suited to
deconvolution problems in the presence of additive Gaussian noise, where the proximity operator associated with the fidelity term takes a closed form \cite{Combettes_PL_2008_journal_proximal_apdmfscvip}.  To minimize a sum of two functions one of which is quadratic, another interesting class of parallel convex optimization algorithms was proposed by Fornasier et al. in \cite{Fornasier_M_2007_j-ip_domain_dmlipsc,Fornasier_M_2009_j-sjna_subspace_cmtvl1m}. 
 However, in a more general context, particularly when the noise is
 not additive Gaussian and a wider class of degradation
   operators is considered, the proximity operator associated with data fidelity
   term does not have a closed form, which prevents the direct use of PPXA and the algorithms in \cite{Fornasier_M_2007_j-ip_domain_dmlipsc}
   and \cite{Fornasier_M_2009_j-sjna_subspace_cmtvl1m}. Therefore,
   other solutions have to be looked for. For Poisson noise, a
   first solution is to resort to the Anscombe transform
   \cite{Dupe_FX_2008_ip_proximal_ifdpniusr}, while a second one  consists of approximating the Poisson data fidelity term with a gradient Lipschitz function \cite{Chaux_C_2009_siims_nested_iafirp}. Both approaches require the use of a nested iterative algorithm \cite{Dupe_FX_2008_ip_proximal_ifdpniusr,Chaux_C_2009_siims_nested_iafirp}, combining forward-backward and Douglas-Rachford iterations. 
Nested algorithms may however appear limited for two main reasons: the
parallelization of the related iterations is difficult, and the number
of functions to be minimized is in practice limited to
three. More recently, approaches related to augmented
  Lagrangian techniques
\cite{Hestenes_M_R_1969_Multipliers_agm,Fortin_M_1983_book_augmented_lmansbvp}
  have been considered in
  \cite{Goldstein_T_2009_siims_split_bml1rp,Setzer_S_2009_Deblurring_pibsbt,Figueiredo_M_2010_t-ip_restoration_piado,He_B_2011_splitting_mscpllc}. These
  methods are well-adapted when the linear operator is a convolution and
  Fourier diagonalization techniques can thus be used, but for more
  general linear degradation operators,  a large-size linear system of equations has to
  be solved numerically at each iteration of the algorithm.

 The
objective of this paper is to propose an adaptation of PPXA to
minimize criteria used in a wide panel of restoration
problems such as those involving a convolution or decimated
  convolution operator using a finite-support kernel and
  non-necessarily additive  Gaussian noise. Decimated convolutions are important in practice since they are
  often encountered in super-resolution problems. To apply proximal methods, it seems that we should be able to compute the proximity operator associated with the fidelity term for a large class of noise distributions. When the proximity operator cannot be easily computed, we will show that
a splitting approach may often be employed to circumvent this difficulty. This is the main contribution of this paper.

Moreover, based on similar splitting techniques and in the spirit of existing works on deconvolution in the presence of Gaussian noise \cite{Bect_J_2004_eccv_unified_vfir,Combettes_PL_2008_journal_proximal_apdmfscvip,BioucasDias_J_2008_Iterative_aiaflipwcr,Wen_YM_2008_sisc_Iterative_abododadfir}, a twofold regularization composed of a sparsity term and a total variation term is performed in order to benefit from each regularization. We will consider this type of hybrid regularization by investigating
different discrete forms of the total variation.

The paper is organized as follows: first, in Section~\ref{sec:context}, we present the considered restoration problem and the general form of the associated criterion to be minimized. Then, in Section~\ref{sec:prox}, the definition and some properties of proximity operators as well as explicit forms related to the data fidelity term in a restoration context and to a discretization of the total variation are provided. 
Section~\ref{sec:algo} introduces an accelerated version of PPXA which allows us to efficiently solve frame-based image recovery problems.
Section~\ref{se:hybrid} shows how the results obtained in the two previous sections can be used for solving restoration problems where
a regularization is performed both in the spatial and in the wavelet domains. 
Finally, in Section~\ref{sec:res}, the effectiveness of the proposed approach 
is demonstrated by experiments for the restoration of images degraded
by a blur (or a decimated blur) with finite-support kernel
  and by a Poisson noise. Some conclusions are drawn in Section~\ref{se:conclu}.

\section{Background}
\label{sec:context}

\subsection{Image restoration}\label{se:restpb}
The degradation model considered throughout this paper is the following:
\begin{equation}
z = \mathcal{D_{\alpha}}(T \overline{y})
\end{equation}
where $\overline{y}$ denotes the original image of global size $N$ degraded by a non-negative valued convolutive operator $T:\RR^N\rightarrow\RR^M$ and contaminated by a noise non necessarily additive, the effect of which is denoted by $\mathcal{D_{\alpha}}$. Here, $\alpha$ is a positive parameter which characterizes the noise intensity. The vector $z\in \RR^M$  represents the observed data of size $M$.
For example, $ \mathcal{D_{\alpha}}$ may denote the addition of a
zero-mean Laplacian noise with standard-deviation $\alpha$,
or the corruption by an independent Poisson noise with scaling parameter  $\alpha$. $T$ represents a convolution or a decimated convolution operator using a finite-support kernel.

Our objective is to recover the image $\overline{y}$ from the observation $z$ 
by using some prior information on its frame coefficients and  
its spatial properties.

\subsection{Frame representation}
In inverse problems, certain physical properties of the target solution
$\overline{y}$ are most suitably expressed in terms of the coefficients
$\overline{x} = (\overline{x}^{(k)})_{1\leq k \leq K}\in \RR^K$
 of its representation
$\overline{y}=\sum_{k=1}^{K}\overline{x}^{(k)}e_k$ with respect to a 
family of vectors $(e_k)_{1\leq k \leq K}$ in the Euclidean space 
$\RR^N$.
Recall that a family of vectors $(e_k)_{1\leq k \leq K}$ in 
$\RR^N$ constitutes a frame if there exist two constants $\underline{\nu}$ and $\overline{\nu}$ 
in $\RPP$ such that \footnote{In finite dimension, the upper bound condition is always satisfied.}
\begin{equation}
\label{e:frame1}
(\forall y\in\RR^N)\qquad\underline{\nu}\|y\|^2\leq\sum_{k=1}^K|\scal{y}{e_k}|^2
\leq\overline{\nu}\|y\|^2.
\;
\end{equation}
The associated frame operator is the injective linear operator $F\colon\RR^N \to \RR^K \colon y\mapsto(\scal{y}{e_k})_{1\leq k \leq K}$,
the adjoint of which is the surjective linear operator $F^{\top}\colon \RR^K\to \RR^N \colon(x^{(k)})_{1\leq k \leq K}\mapsto\sum_{k=1}^K x^{(k)}e_k$.
When $\underline{\nu}=\overline{\nu} = \nu$ in \eqref{e:frame1}, $(e_k)_{1\leq k \leq K}$ is said to be a tight
frame. In this case, we have

\begin{equation}
F^{\top} F = \nu \Idf,
\label{e:tight}
\end{equation} 
where $\Idf$ denotes the identity matrix.
A simple example of a tight frame is the union of $\nu$ orthonormal 
bases, in which case $\underline{\nu}=\overline{\nu}=\nu$. 
For instance, 
a 2D real (resp. complex) dual-tree wavelet decomposition is the 
union of two (resp. four) orthonormal wavelet bases 
\cite{Chaux_C_2006_tip_ima_adtmbwt}. Curvelets \cite{Candes_EJ_2002_as_Recovering_eiipipoocf} constitute another example of tight 
frame. Historically, Gabor frames 
\cite{Daubechies_I_1992_book_ten_lw} have played an important role in many inverse
problems. Under some conditions, contourlets \cite{Do_MN_2005_tip_coutourlet_tctaedmir} also constitute
tight frames.
When $F^{-1}=F^{\top}$, an orthonormal basis is obtained.
Further constructions as well as a detailed account 
of frame theory in Hilbert spaces can be found in \cite{Han_D_2000_book_frames_bgr}.

In such a framework, the observation model becomes
\begin{equation}
z = \mathcal{D_{\alpha}}( T F^{\top} \overline{x} )
\end{equation}
where $\overline{x}$ represents the frame coefficients of the original data ($\overline{y} = F^{\top} \overline{x}\in \RR^N$ is the target data of size $N$). Our objective is now to recover $\overline{x}$ from the observation $z$.

\subsection{Minimization problem}
In the context of inverse problems, the original image can be restored by solving a convex optimization problem of the form:
\begin{equation}
 \mbox{Find}\quad \widehat{x} \in \underset{x \in \RR^K}{\Argmin} \quad \sum_{j=1}^{J}f_j(x)
\label{pb:optinit}
\end{equation}
where $(f_j)_{1\le j\le J}$ are functions of $\Gamma_0(\RR^K)$ (see
\cite{Combettes_PL_2008_journal_proximal_apdmfscvip} and references therein) and the restored image is $\widehat{y} = F^\top\widehat{x}$.

A particular popular case is when $J=2$; the minimization problem thus reduces to the minimization of the sum of two functions which, under a Bayesian framework,
can be interpreted as a fidelity term $f_1$ linked to noise and an a priori term $f_2$ related to some prior probabilistic model put on the frame coefficients
(some examples will be given in Section~\ref{sec:problem}).

In this paper, we are especially interested in the case when $J>2$, which may be fruitful for imposing additional constraints on the target solution.
At the same time, when considering a frame representation (which, as already mentioned, often allows us to better express some properties of the target solution), the convex optimization problem \eqref{pb:optinit} can be re-expressed as:
\begin{equation}
 \mbox{Find}\quad \widehat{x} \in \underset{x \in \RR^K}{\Argmin} \quad \sum_{j=1}^{\JJ}g_j(F^{\top}x) + \sum_{j=\JJ+1}^{J}f_j(x) 
\label{pb:optinit_trame}
\end{equation}
where $(g_j)_{1\le j\le \JJ}$ are functions of $\Gamma_0(\RR^N)$ and $(f_j)_{\JJ+1\le j\le J}$ are functions of $\Gamma_0(\RR^K)$,
related to the image or to the frame coefficients, respectively. The terms for $j\in \{1,\ldots,\JJ\}$ related directly to the pixel values may be the data fidelity term, or a pixel range constraint term, whereas, the functions of indices
$j\in \{\JJ+1,\ldots,J\}$ defined on frame coefficients
are often chosen from some classical prior probabilistic model. For example,
they may correspond to
the minus log-likelihood of independent variables following
generalized Gaussian distributions \cite{Do_MN_2002_tip_GG_wbtruggdakld}.

We will now present convex analysis tools which are useful to deal with such minimization problems.

\section{Proximal tools}

\label{sec:prox}

\subsection{Definition and examples}
A fundamental tool which has been widely employed in the recent convex optimization literature is the proximity operator \cite{Combettes_P_2010_inbook_proximal_smsp} first introduced by Moreau in 1962 \cite{Moreau_JJ_1962_cras_Fonctions_cdeppdueh,Moreau_J_1965_bsmf_Proximite_eddueh}. The proximity operator of $\varphi \in  \Gamma_0(\RR^X)$ is defined as
\begin{equation}
 \prox_\varphi\colon\RR^X \to \RR^X\colon u \mapsto \arg\min_{v \in \RR^X} \Frac12\left\|v-u\right\|^{2} + \varphi(v).
\label{eq:prox}
\end{equation}
Thus, if $C$ is a nonempty closed convex set of $\RR^X$, and $\iota_C$ denotes the indicator function of $C$, \textit{i.e.}, $\forall u\in\RR^X$,  $\iota_C(u)= 0$ if $u\in C$, $+\infty$ otherwise, 
then, $\prox_{\ic}$ reduces to the projection $P_C$ onto $C$. Other examples of proximity operators corresponding to the potential functions of standard log-concave univariate probability densities have been listed in \cite{Combettes_PL_2005_mms_Signal_rbpfbs,Chaux_C_2007_ip_variational_ffbip,Combettes_PL_2008_journal_proximal_apdmfscvip}. Some of them will be used in the paper and we will thus recall the proximity operators of the potentials associated with a Gamma distribution (which is closely related to the Kullback-Leibler divergence \cite{Titterington_D_1987_tmi_iterative_otiisrafe}) and with a generalized Gaussian distribution, before dealing with the Euclidean norm in dimension 2.
\begin{example}\cite{Combettes_PL_2007_istsp_Douglas_rsatncvsr} 
\label{ex:gamd}
Let $\alpha>0$ and set 
\begin{align}
\label{e:gamd}
\varphi\colon&\RR\to\RX \nonumber \\
&\eta\mapsto \begin{cases}
-\chi\ln(\eta)+\alpha\eta,&\text{if}\;\;\chi>0\;\;\text{and}\;\;\eta>0;\\
\alpha\eta,&\text{if}\;\;\chi=0\;\;\text{and}\;\;\eta\ge 0;\\
\pinf, & \text{otherwise}.
\end{cases}
\end{align} 
Then, for every $\eta\in\RR$,
\begin{equation}
\prox_\varphi\eta=\frac{\eta-\alpha+\sqrt{|\eta-\alpha|^2+4\chi}}{2}.
\end{equation}
\end{example}

\begin{example}\cite{Chaux_C_2007_ip_variational_ffbip} 
\label{ex:gg}
Let $\chi>0$, $p\in [1,+\infty[$, and set 
\begin{equation}
\varphi\colon\RR\to\RX\colon\eta\mapsto \chi |\eta|^p. 
\end{equation}
Then, for every $\eta\in\RR$, $\prox_\varphi\eta$ is given by
\begin{equation}
\begin{cases}
\sign(\eta) \max\{|\eta| - \chi,0\} & \mbox{if $p=1$}\\
\eta + \frac{4 \chi}{3 \;.\; 2^{1/3}} \big( (\epsilon - \eta)^{1/3} -
(\epsilon + \eta)^{1/3} \big)\\
\qquad \mbox{where} \quad \epsilon = \sqrt{\eta^2 + 256 \chi^3/729}& \mbox{if $p=\frac{4}{3}$}\\
\eta + \frac{9 \chi^2 \sign(\eta)}{8} \Big(1 - \sqrt{1 + \frac{16
|\eta|}{9\chi^2}}\Big)& \mbox{if $p=\frac{3}{2}$}\\
\frac{\eta}{1+2\chi} & \mbox{if $p=2$}\\
\sign(\eta) \frac{\sqrt{1+12 \chi |\eta|}-1}{6\chi} & \mbox{if $p=3$}\\
\end{cases}
\end{equation}
\end{example}
where $\sign$ denotes the signum function. In Example \ref{ex:gg}, it can be noticed that the proximity operator associated with $p=1$ reduces to a soft thresholding.

\begin{example}\cite{Combettes_PL_2008_journal_proximal_apdmfscvip} 
\label{ex:sqrt}
Let $\mu>0$ and set 
\begin{equation}
\varphi\colon\RR^2\to\RR \colon (\eta_1,\eta_2)\mapsto \mu \sqrt{|\eta_1|^2 + |\eta_2|^2}.
\end{equation}
Then, for every $(\eta_1,\eta_2)\in\RR^2 $,
\begin{equation}
\prox_\varphi(\eta_1,\eta_2)=\begin{cases}
\Big(1-\frac{\mu}{\sqrt{|\eta_1|^2 + |\eta_2|^2}}\Big)(\eta_1,\eta_2),&\text{if}\;\;\sqrt{|\eta_1|^2 + |\eta_2|^2}>\mu;\\
(0,0), & \text{otherwise}.
\end{cases}
\end{equation}
\end{example}

\subsection{Proximity operators involving a linear operator}
\label{ss:proxconv}

We will now study the problem of determining the proximity operator of a function $g=\Psi\circ T$ where $T:\RR^N \rightarrow \RR^M$ is a linear operator,
\begin{equation}\label{eq:defPsi}
\Psi\!\colon \!\RR^M \!\to \!\RX\!\colon\! (u^{(m)})_{1 \le m \le M}\mapsto
\sum_{m=1}^M \psi_m(u^{(m)})
\end{equation}
and, for every $m\in\{1,\ldots,M\}$,
$\psi_m \in \Gamma_0(\RR)$.
As will be shown next, the proximity operator of this function can be 
determined in a closed form for specific cases only. However, $g$ can be decomposed as  a sum of functions for which the proximity operators can be calculated explicitly. 
Firstly, we introduce a property concerning the determination of the proximity operator of the composition of a convex function 
and a linear operator, which constitutes a generalization of \cite[Proposition 11]{Combettes_PL_2007_istsp_Douglas_rsatncvsr} for separable convex functions. The proof of the following proposition is provided in Appendix \ref{ap:proxcomp}. 
\begin{proposition}
\label{p:proxcomp}
Let $X\in \NN^*$, $Y\in \NN^*$, and let $(o_m)_{1\leq m \leq Y}$
be an orthonormal basis of $\RR^{Y}$.
Let $\Upsilon$ be a function such that
\begin{equation}
(\forall u \in \RR^{Y})\qquad
\Upsilon(u) = \sum_{m=1}^Y 
\psi_m(\scal{u}{o_m})
\end{equation}
where $(\psi_m)_{1\leq m \leq Y}$ are functions in $\Gamma_0(\RR)$.
Let $L$ be a matrix in $\RR^{Y\times X}$ such that \begin{equation} \label{eq:diagek}
\underbrace{L L^{\top}}_{D} = \sum_{m=1}^Y \Delta_m o_m o_m^\top
\end{equation}
where $(\Delta_m)_{1\leq m \leq Y}$ is a sequence of positive reals.
\\
Then $\Upsilon\circ L\in\Gamma_0(\RR^X)$ and, for every $v \in \RR^X$
\begin{equation}
\label{e:pfL}
\prox_{\Upsilon\circ L} v= v +L^{\top} D^{-1}\big(\prox_{D\Upsilon}(L v) -L v\big)
\end{equation}
where $D\Upsilon$ is the function
defined by
\begin{equation}
(\forall u \in \RR^{Y})\qquad D\Upsilon(u)=\sum_{m=1}^Y  \Delta_m\psi_m\big(\scal{u}{o_m}\big).
\end{equation}
\end{proposition}

The function $\Psi$ defined in \eqref{eq:defPsi} is separable in the canonical
basis of $\RR^M$.
However, for an abitrary convolutive (or decimated convolutive) operator $L=T$, \eqref{eq:diagek} is generally not satisfied. Nevertheless, assume that $(\II_i)_{1\le i \le I}$
is a partition of $\{1,\ldots,M\}$ in nonempty sets.
For every $i\in \{1,\ldots,I\}$, let $M_i$ be the number of elements in $\II_i$ ($\sum_{i=1}^I M_i = M$) and let
$\Upsilon_i\,:\,\RR^{M_i}\to \RPP\,:\,(u^{(m)})_{m\in \II_i} 
\mapsto \sum_{m\in \II_i} \psi_m(u^{(m)})$.
If, for every $i\in \{1,\ldots,M\}$, $t_i$ is the vector of $\RR^N$ corresponding to the
  $i$-th row vector of $T$, we have then
$g = \sum_{i=1}^I \Upsilon_i \circ T_i $
where $T_i$ is a linear operator from $\RR^N$ to $\RR^{M_i}$ associated with 
a matrix 
\begin{equation}
\begin{bmatrix}
{ t}_{m_1}\; \ldots \; { t}_{m_{M_i}}
\end{bmatrix}^\top
\end{equation}
and $\II_{i} = \{m_1,\ldots,m_{M_i}\}$. The following assumption will play a prominent role in the rest of the paper:

\begin{assumption}\label{a:parti}
For every $i\in \{1,\ldots,I\}$, 
$({ t}_{m})_{m\in\II_i}$ is a family of non zero orthogonal vectors.
\end{assumption}
Then, $g$ can be decomposed as a sum of $I$ functions
$(\Upsilon_i \circ T_i)_{1 \le i \le I}$ where, for every 
$i\in \{1,\ldots,I\}$,  $D_i=T_i T_i^{\top}$
is associated with an invertible 
diagonal matrix $\diag(\Delta_{i,1},\ldots, \Delta_{i,M_i})$.
According to Proposition~\ref{p:proxcomp}, we have then, for every $y\in \RR^N$,
\begin{equation}
\label{eq:proxnewTi}
\prox_{\Upsilon_i \circ T_i} y = y+T_i^{\top}D_i^{-1}\big(\prox_{D_i\Upsilon_i}(T_iy)- T_i y\big).
\end{equation}
\begin{remark}

\begin{enumerate}
 \item Note that Assumption \ref{a:parti} is obviously satisfied when $I=M$, that is 
when, for every $i\in \{1,\ldots,I\}$, $\II_i$ reduces to a singleton.
\item It can be noticed that the application of  $T_i$ or $T_i^{\top}$
  reduces to standard operations in signal processing. For
    example, when $T$ corresponds to a convolutive operator, the application of $T_i$ consists of two steps: a convolution with the impulse response of the degradation filter and a decimation for selected locations ($m\in \mathbb{I}_i$). The application of $T_i^{\top}$ also consists of two steps: an interpolation step (by inserting zeros between data values of indices $m\in \mathbb{I}_i$) followed by a convolution with the filter with conjugate frequency response. 
\end{enumerate}
\end{remark}
The fundamental idea behind the previously introduced partition $(\II_i)_{1\le i
    \le I}$, is to form groups of non-overlapping -- and thus orthogonal -- shifts of the convolution kernel so as to be able to compute the corresponding proximity operators. To reduce the number of proximity operators to be computed, one usually wants to find the smallest integer $I$ such that, for every $i\in \{1,\ldots,I\}$,
$({ t}_{m})_{m\in\II_i}$ is an orthogonal family. For the sake of simplicity, we will consider the case of a 1D deconvolution problem, 
where $N$ represents the original signal size whereas $M$ corresponds
to the degraded signal size, the extension to 2D deconvolution problems being straightforward.
Different configurations concerning the impact of boundary effects on
the convolution operator will be studied: first, we will consider the
case when no boundary effect occurs. Then, boundary effects introduced
by zero padding and by a periodic convolution will be taken into
account. Finally, the special case of decimated convolution
will be considered. $Q$ designates in the sequel the length of the
kernel and $(\theta_{q})_{0\le q<Q}$  its values.

\begin{enumerate}
\item One-dimensional convolutive models without boundary effect.\\
We typically have the following T{\oe}plitz structure:
{\small{\begin{equation}
\begin{bmatrix}
{ t}_{1}^\top\\
\vdots\\
{ t}_{M}^\top
\end{bmatrix}
=
\begin{bmatrix}
\theta_{Q-1} & \ldots 	    & \theta_1 & \theta_0 & 0       & \ldots & 0\\
0            & \ddots       &        &          & \ddots  & \ddots & \vdots\\
\vdots       & \ddots       & \ddots &          &         & \ddots & 0\\
0            & \ldots       & 0      & \theta_{Q-1} & \ldots & \theta_1 & \theta_0
\end{bmatrix}
\end{equation}}}
where $M = N - Q + 1\ge Q$.\\ 
In order to satisfy Assumption~\ref{a:parti}, we can choose $I=Q$ and, for every $i \in \{1,\ldots,I\}$,
\begin{equation}\label{eq:IIi1Dwbe}
 \II_i =\menge{m\in \{1,\ldots,M\}}{(m-i)\!\!\mod I=0}.
\end{equation}
Hence, we have for all $i \in \{1,\ldots,I\}$,
\begin{equation}
\label{eq:diag1Dsimple}
\Delta_{i,1} = \ldots = \Delta_{i,M_i} = \sum_{q=0}^{Q-1}|\theta_q|^2.
\end{equation}
In this case, $g$ can be decomposed as a sum of $Q$ functions,
whose proximity operators can be easily calculated.
\item One-dimensional zero-padded convolutive models.\\
The following T{\oe}plitz matrix is considered:
{\small{\begin{equation}
\begin{bmatrix}
{ t}_{1}^\top\\
\vdots\\
{ t}_{M}^\top
\end{bmatrix} =
\begin{bmatrix}
\theta_0     & 0        & \ldots & \ldots       & \ldots & \ldots & 0\\
\theta_1     & \theta_0   & \ddots &              &        &        & \vdots\\
\vdots       &  \ddots    & \ddots & \ddots       &        &        & \vdots\\
\theta_{Q-1} &        & \ddots & \ddots       & \ddots &        & \vdots\\
0            & \ddots &        &  \ddots      & \ddots & \ddots & \vdots\\
\vdots       & \ddots & \ddots &              &  \ddots & \ddots & 0\\
0            & \ldots & 0      & \theta_{Q-1} & \ldots & \theta_1 & \theta_0
\end{bmatrix}
\end{equation}}}
where $M = N \geq 2Q$. In this case, $I$ can be chosen equal to $Q$ and the index sets
 $(\II_i)_{1\le i \le I}$ are still given by \eqref{eq:IIi1Dwbe}.
However, the diagonal parameters are not all equal as in the previous example.
We have indeed, for every $i \in \{1,\ldots,I\}$,
\begin{equation}
\begin{cases}
\Delta_{i,1} = \sum_{q=0}^{i-1}|\theta_q|^2 \\
\Delta_{i,2} = \ldots = \Delta_{i,M_i} = \sum_{q=0}^{Q-1}|\theta_q|^2.
\end{cases}
\end{equation}

\item One-dimensional periodic convolutive models.\\
In this case, a matrix having a circulant structure \cite{Golub_G_1996_book_matrix_comp} is involved:
{\small{\begin{equation}
\begin{bmatrix}
{ t}_{1}^\top\\
\vdots\\
{ t}_{M}^\top
\end{bmatrix} =
\begin{bmatrix}
  \theta_0   & 0        & \ldots &  0     & \theta_{Q-1} & \ldots & \theta_1\\
\theta_1     & \theta_0   & \ddots &              & \ddots       &  \ddots      & \vdots\\
\vdots       &            & \ddots & \ddots       &        &  \ddots      & \theta_{Q-1}\\
\theta_{Q-1} &        &        & \ddots       & \ddots &        &  0\\
0            & \ddots &        &              & \ddots & \ddots & \vdots\\
\vdots       & \ddots & \ddots &              &        & \ddots & 0\\
0            & \ldots & 0      & \theta_{Q-1} & \ldots & \theta_1 & \theta_0
\end{bmatrix}
\end{equation}}}
where $M = N \ge Q$. 
In order to satisfy Assumption~\ref{a:parti}, we subsequently set 
$I= \min\{i \ge Q\mid (M-i)\!\!\mod Q =0\}$ 
and, for every $i\in \{1,\ldots,I\}$,
\begin{equation}
\II_i = \begin{cases}
\{i\} &\!\!\!\!\!\!\!\!\!\!\!\!\!\!\!\!\!\!\!\!\!\!\!\!\!\!\!\!\!\!\!\!\!\!\!\!\mbox{if $i \leq Q-1$}\\
\menge{m\in \{Q,\ldots,M\}}{(m-i)\!\!\mod Q=0}& \\
 &\!\!\!\!\!\!\!\!\!\!\!\!\!\!\!\!\!\!\!\!\!\!\!\!\!\!\!\!\!\!\!\!\!\!\!\!\mbox{otherwise.}
\end{cases}
\end{equation}

\noindent The diagonal parameters are then given by \eqref{eq:diag1Dsimple}.\\
Another choice which was made in \cite{Pustelnik_N_2009_eusipco_Hybrid_rfdritpopn} is to set 
$I = \min\{i \ge Q \mid M\!\! \mod i = 0\}$ and to proceed as in \eqref{eq:IIi1Dwbe} and \eqref{eq:diag1Dsimple}. This solution may be preferred due to its simplicity, when the resulting value of $I$ is small.

\item One-dimensional $d$-decimated zero-padded convolutive models.\\
We get the following matrix of $\RR^{M\times N}$ where $N = Md \ge 2Q$.
\end{enumerate}
{\small{\begin{equation}
\begin{bmatrix}
{ t}_{1}^\top\\
\vdots\\
{ t}_{M}^\top
\end{bmatrix} =
\begin{bmatrix}
\mbox{\begin{tabular}{c c c c c c c c c c c}
$\theta_{d-1}$ & $\ldots$ & $\theta_0$   &  $0$   &  $\ldots$		&  $\ldots$    & $\ldots$  &$\ldots$ 	& 	$\ldots$ 	 & $\ldots$ & $0$\\ 
$\theta_{2d-1}$ & $\ldots$  & $\theta_{d-1}$ & $\ldots$ &  $\theta_0$ & $0$          & 	 $\ldots$	&  $\ldots$      & 	$\ldots$  & $\ldots$ & $0$\\
 $\vdots$  & $\ddots$ &  &  $\ddots$ &  & $\ddots$ &$\ddots$  &  &  & &  $\vdots$\\
 &  & $\ddots$ &   & $\ddots$ &  & $\ddots$ & $\ddots$ &  & &  $\vdots$\\
 0 &  &  & $\ddots$  &  & $\ddots$ &  & $\ddots$ & $\ddots$ & &  $\vdots$\\
  $\vdots$ & $\ddots$ &  &   & $\ddots$ &  & $\ddots$ &  & $\ddots$ & $\ddots$&  $\vdots$\\
 $\vdots$ &  & $\ddots$ &   &  &$\ddots$  &  & $\ddots$ &  & $\ddots$ & 0\\
 $0$       & $\ldots$   &$\ldots$      & $0$      	& $\theta_{Q-1}$& $\ldots$ 	& $\theta_{2d-1}$&  $\ldots$	& $\theta_{d-1}$&  $\ldots$ &  $\theta_0$
\end{tabular}}
\label{eq:matlong}
 \end{bmatrix}.
\end{equation}}}

\begin{itemize}
\item[]
In order to satisfy Assumption~\ref{a:parti}, we subsequently set 
$I= \big \lceil \frac{Q}{d}\big\rceil$ and the index sets
 $(\II_i)_{1\le i \le I}$ are still given by \eqref{eq:IIi1Dwbe}. 
We have indeed, for every $i \in \{1,\ldots,I\}$,
\begin{equation}
\label{eq:decconv_diag}\begin{cases}
\Delta_{i,1} = \sum_{q=0}^{\min(i d,Q)-1}|\theta_q|^2 \\ 
\Delta_{i,2} = \ldots = \Delta_{i,M_i} = \sum_{q=0}^{Q-1}|\theta_q|^2.                        
                       \end{cases}
\end{equation}
Note that, when $d\geq Q$, $({ t}_{m})_{m\in\{1,\ldots,M\}}$ is an
orthogonal family, and thus $I=1$.
\end{itemize}

\begin{remark}
In the previous example (the non-decimated example being a special
case when $d=1$), the computational complexity of applying
each operator $T_i$ or $T_i^{\top}$ with $i\in \{1,\ldots,I\}$ is $O(Md)$
and we have about $Q/d$ proximity operators $\prox_{\Upsilon_i\circ T_i}$
to compute. Assuming a complexity $O(M_i)$ for computing 
$\prox_{D_i\Upsilon_i}$, the overall computational complexity
is $O(M(2Q+1))$. 
In turn, if we choose $I=M$, the complexity of computation of $T_i$ or $T_i^{\top}$
is $O(Q)$, but we have about $M$ proximity operators $\prox_{\Upsilon_i\circ T_i}$
to compute. Thus, the overall computational complexity remains of the same order as previously.
This means that limiting the number of proximity operators
to be computed has no clear advantage in terms of computational complexity, but
it allows us to reduce the memory requirement (gain of a factor $Md/Q$ for
the storage of the results of the proximity operators).
\end{remark}

\subsection{Discrete forms of total variation and associated proximity operator}
\label{ss:proxtv}
Total variation \cite{Rudin_L_1992_tv_atvmaopiip}  represents a powerful regularity measure in image restoration for recovering piecewise homogeneous
areas with sharp edges \cite{Rudin_L_1994_tv_birflc,Malgouyres_F_2002_jip_math_amctvwir,Aujol_JF_2006_ijcv_structure_tidmaps,Weiss_P_2009_sjsc_efficient_stvmcip,Bresson_X_2007_j-ipi_fast_dmv}. Different versions of discretized total variation can be found in the literature \cite{Rudin_L_1992_tv_atvmaopiip,Chambolle_A_2004_jmiv_TV_aaftvmaa,Combettes_PL_2008_journal_proximal_apdmfscvip}. Our objective here
is to consider discrete versions for which the proximity
operators can be easily computed. The main idea will be to split the
total variation term in a sum of functions the proximity
  operators of which have a closed form.
The considered form of the total variation of a digital image
$y = (y_{n_1,n_2})_{0 \le n_1 < N_1,0\le n_2 < N_2}\in \RR^{N_1\times N_2}$ is
\begin{equation}
\label{eq:tv_general}
\tv(y) 
=\sum_{n_1=0}^{N_1-P_1}\sum_{n_2=0}^{N_2-P_2} \rho_{\rm{tv}}\big((y_{h})_{n_1,n_2},(y_{v})_{n_1,n_2}\big),
\end{equation}
where $\rho_{\rm{tv}}\in  \Gamma_0(\RR^2)$, and $y_{h}$ and $y_{v}$ are two discrete gradients computed in orthogonal directions through FIR filters with impulse responses of size $P_1\times P_2$. More precisely, in the above expression,
we have
\begin{equation}
\begin{cases}
(y_{h})_{n_1,n_2} & = \trace(\hh^\top Y_{n_1,n_2})\\
(y_{v})_{n_1,n_2} & = \trace(\vv^\top Y_{n_1,n_2}) 
\end{cases}
\end{equation}
where $\hh\in\RR^{P_1\times P_2}$ and $\vv\in\RR^{P_1\times P_2}$
are the filter kernel matrices here assumed to have unit Frobenius norm, and for every $(n_1,n_2) \in \{0,\ldots,N_1-P_1\} \times
\{0,\ldots,N_2-P_2\}$, 
$Y_{n_1,n_2} = (y_{n_1+p_1,n_2+p_2})_{0 \le p_1 < P_1, 0 \le p_2 < P_2}$ denotes a block of $P_1\times P_2$ neighbouring pixels.
Since the proximity operator associated with the so-defined total variation
does not take a simple expression in general,  \eqref{eq:tv_general}
can be split in ``block terms'' by following an approach similar to
that in
Section~\ref{ss:proxconv}:
\begin{equation}
\label{eq:tv_particulier}
(\forall y \in \RR^{N_1\times N_2}) \qquad \tv(y) = \sum_{p_1=0}^{P_1-1}\sum_{p_2=0}^{P_2-1}\tv_{p_1,p_2}(y)
\end{equation}
where, for every $p_1 \in \{0,\ldots, P_1-1\}$
and $p_2 \in \{0,\ldots, P_2-1\}$,
\begin{equation}
\tv_{p_1,p_2}(y) =  \!\!\!\!\!\!\sum_{n_1=0}^{\lfloor\frac{N_1-p_1}{P_1}\rfloor -1}\sum_{n_2=0}^{\lfloor\frac{N_2-p_2}{P_2}\rfloor-1}  \!\!\!\!\!\!\rho_{\rm{tv}}\big((y_{h})_{n_1,n_2}^{p_1,p_2},(y_{v})_{n_1,n_2}^{p_1,p_2}\big)
\end{equation}
and the notation $(\cdot)_{n_1,n_2}^{p_1,p_2} = (\cdot)_{P_1
  n_1+p_1,P_2 n_2 +p_2}$ has been used. A closed form expression
for the proximity operator of the latter function can be
derived as shown below (the proof is provided in Appendix~\ref{ap:proxtv}).

\begin{proposition}
\label{p:proxtv}
Under the assumption that $\trace(\hh \vv^\top) = 0$, 
for every 
\begin{equation*}
y = (y_{n_1,n_2})_{0 \le n_1 < N_1,0\le n_2 < N_2}
\in \RR^{N_1\times N_2} 
\end{equation*}
and $\mu > 0$, we have
\begin{equation}
(\forall (p_1,p_2) \in \{0,\ldots,P_1-1\}\times
\{0,\ldots,P_2-1\}) \quad \prox_{\mu\tv_{p_1,p_2}}y  = (\pi_{n_1,n_2})_{0 \le n_1 < N_1,0\le n_2 < N_2}
\end{equation}
where, for every $(n_1,n_2) \in \{0,\ldots,\lfloor\frac{N_1-p_1}{P_1}\rfloor-1\}
\times \{0,\ldots,\lfloor\frac{N_2-p_2}{P_2}\rfloor-1\}$,
\begin{equation}
\label{eq:proxtv1}
(\pi_{P_1 n_1+p_1+p_1',P_2n_2+p_2+p_2'})_{0\le p_1'<P_1,0\le p_2'<P_2}=
(\beta_{n_1,n_2}^{p_1,p_2}-h_{n_1,n_2}^{p_1,p_2}) \hh + (\kappa_{n_1,n_2}^{p_1,p_2}-v_{n_1,n_2}^{p_1,p_2})\vv + Y_{n_1,n_2}^{p_1,p_2}
\end{equation}
with
\begin{align}
&h_{n_1,n_2}^{p_1,p_2} = \trace(\hh^\top Y_{n_1,n_2}^{p_1,p_2}),\quad
v_{n_1,n_2}^{p_1,p_2} = \trace(\vv^\top Y_{n_1,n_2}^{p_1,p_2}) \label{eq:yv}\\
&(\beta_{n_1,n_2}^{p_1,p_2},\kappa_{n_1,n_2}^{p_1,p_2}) =  
\prox_{\mu \;\rho_{\rm tv}}(h_{n_1,n_2}^{p_1,p_2},v_{n_1,n_2}^{p_1,p_2}),\label{eq:proxtv2}
\end{align}
and 
$(\forall (n_1,n_2) \in \{0,\ldots,N_1-1\} \times
\{0,\ldots,N_2-1\})$
\begin{equation}\label{eq:pibord}
\pi_{n_1,n_2} = y_{n_1,n_2} \;\;\mbox{if}\;\;\begin{cases}
 \mbox{$n_1 <p_1$ or}\\
\mbox{$n_2 < p_2$ or}\\
\mbox{$n_1 \ge P_1\lfloor\frac{N_1-p_1}{P_1}\rfloor$ or}\\
\mbox{$n_2 \ge P_2\lfloor\frac{N_2-p_2}{P_2}\rfloor$}.                             
                            \end{cases}
\end{equation}
\end{proposition}

The result in Proposition~\ref{p:proxtv} basically means
  that, for a given value of $(p_1,p_2)\in \{0,\ldots,P_1-1\}\times
  \{0,\ldots,P_2-1\}$, the image is decomposed into
non-overlapping blocks $Y_{n_1,n_2}^{p_1,p_2}=$\linebreak $(y_{P_1
  n_1+p_1+p_1',P_2n_2+p_2+p_2'})_{0\le p_1'<P_1,0\le p_2'<P_2}$ of $P_1\times P_2$ pixels. Eq. \eqref{eq:proxtv1} then provides the expression of the proximity
operator associated with each one of these blocks, whereas
\eqref{eq:pibord} deals with boundary effects.

\begin{remark}
The above result offers some degrees of freedom in the definition of the discretized total variation for the choices
of the function $\rho_{\rm tv}$ and of the gradient filters.
\begin{itemize}
\item Two classical choices for the function 
$\rho_{\rm tv}$ \cite{Rudin_L_1992_tv_atvmaopiip} are the following:
\begin{enumerate}
 \item If $\rho_{\rm tv}\colon (\eta_1,\eta_2)\mapsto|\eta_1|+|\eta_2|$ then, an anisotropic
form is obtained. According to Example \ref{ex:gg}, 
\eqref{eq:proxtv2} reduces to
\begin{equation}
\begin{cases}
\beta_{n_1,n_2}^{p_1,p_2} = \sign(h_{n_1,n_2}^{p_1,p_2} ) \max(|h_{n_1,n_2}^{p_1,p_2} |-\mu,0)\\
\kappa_{n_1,n_2}^{p_1,p_2} = \sign(v_{n_1,n_2}^{p_1,p_2} ) \max(|v_{n_1,n_2}^{p_1,p_2} |-\mu,0).
\end{cases}
\end{equation}
 \item If $\rho_{\rm tv}\colon (\eta_1,\eta_2)\mapsto\sqrt{(\eta_1)^2+(\eta_2)^2}$, then 
the standard isotropic form is found. The proximity operator
involved in \eqref{eq:proxtv2} is given in Example \ref{ex:sqrt}.
\end{enumerate}
\item Some examples of kernel matrices $\hh$ and $\vv$ satisfying the assumptions
of Proposition \ref{p:proxtv} are as follows:
\begin{enumerate}
 \item Roberts filters such that
$\hh =$ {\small{$ \begin{bmatrix}
           -1/\sqrt{2} & 0 \\ 0 & 1/\sqrt{2}
         \end{bmatrix}$}}
and
 $\vv=$  {\small{$\begin{bmatrix}
          0 &-1/\sqrt{2} \\  1/\sqrt{2} & 0
         \end{bmatrix}$}}
were investigated in \cite{Combettes_PL_2008_journal_proximal_apdmfscvip}.
\item Finite difference filters can be used, which are such that
$\hh = \vv^{\top} =   {\small{\begin{bmatrix}
0 & 0 & 0\\
            -1/\sqrt{2} & 0 &1/\sqrt{2}\\
0 & 0 & 0\\
         \end{bmatrix}}}$.
\item Prewitt filters also satisfy the required assumptions. They are defined by\\
$\hh = \vv^{\top} =  {\small{ \begin{bmatrix}
             -1/\sqrt{6}& 0& 1/\sqrt{6} \\ -1/\sqrt{6}& 0& 1/\sqrt{6} \\ -1/\sqrt{6} & 0& 1/\sqrt{6}
         \end{bmatrix}}}$.
\item Sobel filters such that\\
  $\hh = \vv^{\top} = {\small{\begin{bmatrix}
             -1/\sqrt{12}& 0& 1/\sqrt{12} \\ -2/\sqrt{12}& 0& 2/\sqrt{12} \\ -1/\sqrt{12} & 0& 1/\sqrt{12}
         \end{bmatrix}}}$
are possible choices too.
\end{enumerate}
\end{itemize}
\end{remark}

\section{Proposed algorithm}\label{sec:algo}
In the class of convex optimization methods, an algorithm recently proposed
in \cite{Combettes_PL_2008_journal_proximal_apdmfscvip} appears 
well-suited to solve the class of the minimization problems formulated as in 
Problem~\eqref{pb:optinit}. However, when synthesis frame
representations are considered (Problem~\eqref{pb:optinit_trame}) and
when the function number $S$ is large, the frame analysis and
synthesis operators have to be applied several times in the algorithm
which induces a long computation time. In this section,
we briefly recall the Parallel ProXimal Algorithm and its
convergence properties. Then, we propose an improved version of PPXA
to efficiently solve Problem~\eqref{pb:optinit_trame}.

\subsection{Parallel ProXimal Algorithm (PPXA)}
An equivalent formulation of the convex optimization problem~\eqref{pb:optinit} is:
\begin{equation}
 \mbox{Find}\quad \widehat{x} \in  \underset{\substack {x_1 \in \RR^K,\ldots, x_J \in \RR^K\\x = x_1=\ldots = x_J}}{\Argmin} \quad \sum_{j=1}^{J}f_j(x_j).
\label{pb:optinitgen}
\end{equation}
This formulation was used in \cite{Combettes_PL_2008_journal_proximal_apdmfscvip} to derive Algorithm~\ref{algo:PPA}.

\begin{algorithm}[!ht]
\caption{General form of PPXA \label{algo:PPA}}
\begin{algorithmic}
\STATE Set $\gamma \in \RPP$.
\STATE For every $j\in\{1,\ldots,J\}$, set $(\omega_j)_{1\leq j\leq J}\in]0,1]^J$ such that $\sum_{j=1}^J \omega_j =1$.
\STATE Set $(u_{j,0})_{1\leq j\leq J} \in (\RR^K)^{J}$ and $x_0 = \; \sum_{j=1}^{J} \omega_j u_{j,0}$.
\STATE {$\mbox{For} \; \ell=0,1,\ldots$}
\STATE 
$\left \lfloor \begin{array}{l}
\mbox{For} \;  j=1,\ldots,J \\
\lfloor \quad p_{j,\ell} = \prox_{\gamma f_j/\omega_j}  u_{j,\ell} +  a_{j,\ell} \\
p_\ell = \sum_{j=1}^{J}\omega_j p_{j,\ell} \\
\mbox{Set } \lambda_\ell\in \;] 0,2[ \\
\mbox{For} \; j=1,\ldots,J \\
\lfloor \quad u_{j,\ell+1} = u_{j,\ell} +\lambda_\ell \; (2\; p_\ell - x_\ell -p_{j,\ell}) \\
x_{\ell+1} = x_\ell + \lambda_\ell  \;(p_\ell -x_\ell)	
\end{array} \right.$
\end{algorithmic}
\end{algorithm}

PPXA involves real constants $\gamma$ and
$(\omega_j)_{1\le j \le J}$, and, at each iteration $\ell \in \NN$,
a relaxation parameter $\lambda_\ell$. It also includes
possible error terms
$(a_{j,\ell})_{1\le j \le J}$ in the computation of the
proximity operators, which shows the numerical stability of the algorithm.
The sequence $(x_{\ell})_{\ell\geq1}$ generated by Algorithm
  \ref{algo:PPA} can be shown to converge to a solution to
Problem~\eqref{pb:optinitgen} (or equivalently to Problem~\eqref{pb:optinit}) under the following assumption \cite{Combettes_PL_2008_journal_proximal_apdmfscvip}.

\begin{assumption}\ 
\begin{enumerate}
 \item \label{a:1} $\lim_{\Vert x\Vert\rightarrow+\infty} f_1(x)+\ldots +f_J(x) = +\infty$.
 \item \label{a:2} $\cap_{j=1}^{J}\rint \dom f_j \neq \emp$.\footnote{The relative interior of a set $S$ of $\RR^X$ is designated by $\rint S$ and the domain of a function $f: \RR^X\rightarrow]-\infty,+\infty]$ is $\dom f = \{x\in \RR^X |f(x)<+\infty\}$.}
\item  \label{a:3} $(\forall j\in\{1,\ldots,J\}) \; \sum_{\ell\in\NN}\lambda_\ell \;\Vert a_{j,\ell}\Vert<+\infty$.
\item  \label{a:4} $\sum_{\ell\in\NN}\lambda_\ell \;(2-\lambda_\ell)= \pinf$.
\end{enumerate}
\label{ass:conv_spingarn}
\end{assumption}

\begin{remark}
The fact that the algorithm involves several parameters should not be
viewed as a weakness since the convergence is guaranteed for any
choice of these parameters under the previous assumption.
These parameters bring out flexibility in  PPXA in the sense that an appropriate
choice of them (typical values will be indicated in Section~\ref{sec:res}) 
may be beneficial to the convergence speed.
\end{remark}

Consider now Problem \eqref{pb:optinit_trame} where a tight frame is
employed
($F^\top F = \nu \,\Idf$).
By setting $(\forall j \in \{1,\ldots,S\})$ $f_j = g_j \circ F^{\top}$ and by invoking Proposition \ref{p:proxcomp} with $L=F^{\top}$ and $D = \nu \,\Idf$, the iterations of Algorithm~\ref{algo:PPA} become as described in Algorithm~\ref{algo:PPA_frame}.

\begin{algorithm}[!ht]
\caption{PPXA iterations for Problem \eqref{pb:optinit_trame} \label{algo:PPA_frame}}
\begin{algorithmic}
\STATE {$ \mbox{For} \; \ell=0,1,\ldots$}
	\STATE 
$\left \lfloor \begin{array}{l}
\mbox{For} \; j=1,\ldots,\JJ\\
\lfloor  \quad p_{j,\ell} = u_{j,\ell} + \\ \qquad\qquad \frac{1}{\nu} F \big(\prox_{\nu\gamma g_j/\omega_j}  (F^{\top}u_{j,\ell}) - F^{\top}u_{j,\ell}\big) +  a_{j,\ell}\\
\mbox{For} \;j=\JJ + 1,\ldots,J\\
\lfloor \quad p_{j,\ell} = \prox_{\gamma f_j/\omega_j}  u_{j,\ell} +  a_{j,\ell}\\
p_\ell = \sum_{j=1}^{J}\omega_j p_{j,\ell}\\
\mbox{Set } \lambda_\ell\in \;] 0,2[\\
\mbox{For} \;j=1,\ldots,J\\
\lfloor  \quad u_{j,\ell+1} = u_{j,\ell} +\lambda_\ell\; (2\; p_\ell - x_\ell -p_{j,\ell})\\
x_{\ell+1} = x_\ell + \lambda_\ell(p_\ell -x_\ell)
 \end{array} \right.$
\end{algorithmic}
\end{algorithm}

However, the first loop can be costly in terms of computational complexity because it requires to apply $\JJ$ times the operators $F$ and $F^{\top}$ at each iteration. We will now see
how it is possible to speed up these iterations.

\subsection{Accelerated version of PPXA}
\label{ss:accPPXA}
In Algorithm \ref{algo:accPPA}, we propose  an adaptation of PPXA in order to reduce its computational load
by limiting the number of times the operators $F$ and $F^{\top}$ are applied. Details concerning the derivation of this algorithm can be found in Appendix \ref{algop:PPA_acc}.

\begin{algorithm}[!ht]
\caption{Accelerated PPXA \label{algo:accPPA}}
\begin{algorithmic}
\STATE Let $\gamma \in \RPP$.
\STATE For every $j\in\{1,\ldots,J\}$, set $(\omega_j)_{1\leq j\leq J}\in]0,1]^J$ such that $\sum_{j=1}^J \omega_j =1$.
\STATE Set $(u_{j,0})_{1\leq j\leq J} \in (\RR^K)^{J}$ and $x_0 = \sum_{j=1}^{J} \omega_j u_{j,0}$.
\STATE For every $j\in\{1,\ldots,\JJ\}$, set $v_{j,0} = F^\top u_{j,0}$ and $u_{j,0}^\perp = u_{j,0}-\frac{1}{\nu}
F v_{j,0}$.
\STATE {For $\ell =0,1,\ldots$}
	 \STATE $\left \lfloor \begin{array}{l} \mbox{For} \; j=1,\ldots,\JJ\\
	 	\lfloor \quad q_{j,\ell} =  \frac{1}{\nu} \prox_{\nu \gamma g_j/\omega_j } v_{j,\ell}+\widetilde{a}_{j,\ell}\\
	 	\mbox{For} \; j=\JJ+1,\ldots,J \\
		\lfloor \quad p_{j,\ell} = \prox_{\gamma f_j/\omega_j}  u_{j,\ell}+a_{j,\ell} \\
	 	p_\ell = \sum_{j=1}^{\JJ}\omega_j {u}_{j,\ell}^{\perp} \\ \qquad \qquad + F \sum_{j=1}^{\JJ}\omega_j q_{j,\ell} + \sum_{j=\JJ+1}^{J}\omega_j p_{j,\ell}\\
	 	r_\ell = 2\;p_\ell - x_\ell ;\quad
		\widetilde{r}_\ell = F^{\top} r_\ell ;\quad
	 	{r}_{\ell}^{\perp} = r_\ell - \frac{1}{\nu}F\widetilde{r}_\ell\\
	 	\mbox{Set } \lambda_\ell\in \;] 0,2[\\
	 	\mbox{For} \;j=1,\ldots,\JJ\\
	 	\left\lfloor
		\begin{tabular}{c}
		${u}_{j,\ell+1}^{\perp} = {u}_{j,\ell}^{\perp} +\lambda_\ell \;({r}_{\ell}^{\perp}  -{u}_{j,\ell}^{\perp})$\\
		$v_{j,\ell+1} = v_{j,\ell} +\lambda_\ell \; (\widetilde{r}_\ell - \nu q_{j,\ell})$
		\end{tabular}
		\right. \\
		\mbox{For} \;j=\JJ+1,\ldots,J \\
		\lfloor \quad u_{j,\ell+1} = u_{j,\ell} +\lambda_\ell \; (r_\ell -p_{j,\ell})\\
	 	x_{\ell+1} = x_\ell + \lambda_\ell (p_\ell -x_\ell)\\
		\end{array} \right. $
\end{algorithmic}
\end{algorithm}

Let us make the following assumption:
\begin{assumption}\ 
\label{as:errorterms}
\begin{enumerate}
\item \label{a:1acc} $\lim_{\Vert x\Vert\rightarrow\pinf} g_1(F^{\top}x)+\ldots+g_\JJ(F^{\top}x)+ f_{\JJ+1}(x)+\ldots +f_J(x) = \pinf$.
\item  \label{a:2acc} $\big(\cap_{j=1}^{\JJ}\rint \dom (g_j\circ F^{\top})\big)\bigcap \big(\cap_{j=\JJ+1}^{J}\rint \dom f_j\big)\neq \emp$.
 \item \label{a:3acc}
$(\forall j\in\{1,\ldots,\JJ\}) \; \sum_{\ell\in\NN}\lambda_\ell \;\Vert \widetilde{a}_{j,\ell}\Vert<+\infty$ and
$(\forall j\in\{\JJ+1,\ldots,J\}) \; \sum_{\ell\in\NN}\lambda_\ell \;\Vert a_{j,\ell}\Vert<+\infty$.
\item  \label{a:4acc} $\sum_{\ell\in\NN}\lambda_\ell \;(2-\lambda_\ell)= \pinf$.
\end{enumerate}
\end{assumption}
Then,
Algorithm \ref{algo:accPPA} converges to a solution to Problem \eqref{pb:optinit_trame}. In addition, this algorithm requires only 3 applications of $F$
or $F^{\top}$ at each iteration. 
Hence, a gain w.r.t. Algorithm \ref{algo:PPA_frame} is obtained as soon as $S \ge 2$. This fact will be illustrated by our simulation results in Section~\ref{se:simspeed}.

\section{Application to restoration}
\label{sec:problem}

\subsection{Hybrid regularization} \label{se:hybrid}
In restoration problems, one of the terms in the criterion to be minimized
usually is a fidelity term measuring some distance between 
the image degraded by the operator $T$
and the observed data $z$. We will assume that this function takes
the form $g = \Psi\circ T$ where $\Psi\in\Gamma_0(\RR^M)$.
In the case of data corrupted by a additive zero-mean white Gaussian noise with variance $\alpha$, a standard choice for 
$\Psi$ is a quadratic function such that $\Psi  = \frac{1}{2\alpha}\Vert \cdot -z\Vert^2$.
Then, the associated proximity operator of $g$ can be computed explicitly (see \cite{Combettes_PL_2008_journal_proximal_apdmfscvip}).
In the case of data contaminated by an independent Poisson noise with scaling parameter $\alpha$, a standard choice is $\Psi = D_{\rm KL}(z,\alpha \;\cdot)$
where $D_{\rm KL}$ is the generalized Kullback-Leibler divergence \cite{Titterington_D_1987_tmi_iterative_otiisrafe, Fessler_JA_1995_tip_Hybrid_ppofftirfts,Zheng_J_2000_j-tip_parallelizable_btargc,Combettes_PL_2007_istsp_Douglas_rsatncvsr, Figueiredo_M_2010_t-ip_restoration_piado, Setzer_S_2009_Deblurring_pibsbt,Chouzenoux09} such that,
\begin{equation}
\label{eq:psi}
\big(\forall u=(u^{(m)})_{1\leq m \leq M}\in \RR^M \big), \quad \Psi(u) = \sum_{m=1}^{M} \psi_m(u^{(m)})
\end{equation}
and
\begin{equation}
\psi_m(u^{(m)}) =
\begin{cases}
\displaystyle \alpha u^{(m)}-z^{(m)}
 + z^{(m)} &\ln \Big(\frac{z^{(m)}}{\alpha u^{(m)}}\Big)
\\
& \!\!\!\mbox{if $z^{(m)}>0$ and $u^{(m)}>0$,}\\
\displaystyle \alpha u^{(m)}
& \! \!\!\mbox{if $z^{(m)}=0$ and $u^{(m)}\geq 0$,}\\
+\infty & \! \!\!\mbox{otherwise.}
\end{cases}
\label{eq:KL}
\end{equation}
The proximity operator of $\Psi$ can then be derived from Example \ref{ex:gamd}.

Concerning regularization functions, a standard choice of penalty
function in the wavelet domain is:   $(\forall x = (x^{(k)})_{1\leq k\leq K}\in\RR^K)$, $\Phi(x) = \sum_{k =1}^{K}\phi_k(x^{(k)})$ where, for every $k\in \{1,\ldots,K\}$, $\phi_k$ is a finite function of $\Gamma_0(\RR)$ such that $\lim_{|x^{(k)}|\to \pinf} \phi_k(x^{(k)}) = \pinf$. Power functions as
in Example~\ref{ex:gg} are often chosen for $(\phi_k)_{1\le k \le K}$
(see e.g. \cite{Daubechies_I_2004_cpamath_iterative_talipsc,
Chaux_C_2007_ip_variational_ffbip}).
The main problem with wavelet regularization is the occurence of some visual artefacts 
(e.g. ringing artefacts), some of which can be reduced by increasing the redundancy of the representation. Another popular type of regularization that can be envisaged consists of employing a total variation measure \cite{Rudin_L_1992_tv_atvmaopiip}. Its major drawback is the generation of staircase-like effects in the recovered images.
To combine the advantages of both regularizations, we propose
to:
 \begin{equation}
\label{eq:minimization}
\mbox{Find}\quad \widehat{x} \in  \underset{x \in \RR^K}{\Argmin} \quad \Psi(T F^{\top}x)  + \mu\, \tv(F^{\top} x) +\iota_C(F^{\top} x)+ \vartheta \; \Phi(x).
\end{equation}
As already mentioned, $\Phi$ corresponds to the regularization term operating in the wavelet domain. $\tv$ represents a discrete total variation term 
as defined by \eqref{eq:tv_general}. 
Finally, $\iota_C$ is the indicator function of a nonempty closed convex set $C$ of $\RR^N$ (for example, related to support or value range contraints).
This kind of objective function was also recently investigated in \cite{Combettes_PL_2008_journal_proximal_apdmfscvip} but the approach was restricted to the use of a quadratic data fidelity term and of a specific form of the total variation
term.

The non-negative real parameters $\vartheta$ and $\mu$ control the degree of smoothness in the wavelet and in the space domains, respectively.

The main difficulty in applying Algorithm \ref{algo:PPA} to our restoration problem is that it requires to compute the proximity operators associated with each of the four terms in \eqref{eq:minimization}. In general, closed forms of the
proximity operators are known only for the indicator function $\iota_C$ and for $\Phi$ \cite{Chaux_C_2007_ip_variational_ffbip}. However, as explained in
Section \ref{ss:proxconv},
provided that the function $\Psi$ is separable, the data fidelity term can be decomposed as a sum of $I$ functions $(\Upsilon_i\circ T_i)_{1\le i \le I}$ for which the proximity operators can be calculated according to \eqref{eq:proxnewTi}. Similarly, by using the results in Section~\ref{ss:proxtv}, the $\tv$
function can be split in $P_1P_2$ functions $(\tv_{p_1,p_2})_{0\le p_1 < P_1,0\le p_2 < P_2}$, the proximity
operators of which are given by Proposition \ref{p:proxtv}. 
Algorithm~\ref{algo:accPPA} can then be applied with $\JJ = I+P_1P_2+1$
and $J = I+P_1P_2+2$. In the present case, it can be noticed that if $\vartheta > 0$,
Assumption \ref{as:errorterms}~\ref{a:1}) is satisfied.
In addition,  Assumption \ref{as:errorterms}~\ref{a:2}) is fulfilled
$\Big(\cap_{i=1}^I \menge{y \in \RR^N}{T_i y \in \rint \dom
    \Upsilon_i}\Big) \bigcap \rint C \neq \emp$ 
 (since $\dom \Phi = \RR^K$ and $(\forall(p_1,p_2)\in \{0,\ldots,P_1-1\}\times \{0,\ldots,P_2-1\})$ $\dom \tv_{p_1,p_2} = \RR^N$).
This condition is verified if $]0,\pinf[^M\subset\dom\Psi$ and $C =
[0,255]^N$ since for every $i\in\{1,\ldots ,I\}$, $T_i$ has been
assumed non-negative real valued in Section \ref{se:restpb}, and with non-zero lines (see Assumption~\ref{a:parti}).

\subsection{Experimental results for convolved data in the presence of Poisson noise}
\label{sec:res}
In our simulations, we will be first interested in studying the performance in terms of convergence rate of the accelerated version of PPXA. Algorithms \ref{algo:PPA_frame} and \ref{algo:accPPA} are implemented by setting $\gamma = 50$,
$\lambda_{\ell} \equiv 1.6$ and, for every $j \in \{1,\ldots,J\}$,
\begin{equation}
\omega_j = \begin{cases}
\frac{1}{4I} & \mbox{if $1 \le j \le I$}\\
\frac{1}{4P_1P_2} & \mbox{if $I+1 \le j \le I+P_1P_2$}\\
\frac{1}{4} & \mbox{otherwise}
\end{cases}
\end{equation}
if $(g_j)_{1 \le j \le I}$ are the functions corresponding 
to the decomposition
of the data fidelity term and\linebreak $(g_j)_{I+1 \le j \le I+P_1P_2}$ correspond
to the decomposition of $\tv$. 
The weights are thus chosen to provide equal
contributions to the four functions in Criterion (42). 
For the first and second functions which are splitted, the corresponding 1/4
weight is further subdivided in a uniform manner.
Note however that the behaviour of the algorithm did not appear to be very
sensitive to an accurate choice of these parameters. 
A comparison between the different total variation regularization terms defined in Section~\ref{ss:accPPXA} will also be made. 
Another discussion will be held concerning the boundary effects. Two
cases will be considered: the use of a periodic convolution and then,
of a convolution with zero-padding. 
Results for a decimated
  convolution will also be presented. 
Finally, the interest in
combining total variation and wavelet regularization terms will be
shown with respect to classical regularizations. 
A tight frame version of the dual-tree transform (DTT) proposed in \cite{Chaux_C_2006_tip_ima_adtmbwt} ($\nu = 2$) using Symlets of length 6 over 3 resolution levels is employed.
We choose potential functions of the form: for every $k\in \{1,\ldots,K\}$, $\phi_k = \chi_k | \cdot |^{p_k}$ where $\chi_k> 0$ and $p_k \in \{1,4/3,3/2,2\}$,
the proximity operators of which are given by Example \ref{ex:gg}.

\subsubsection{Convergence rate comparison between PPXA and its accelerated version}
\label{se:simspeed}
\begin{figure*}[htbp] 
\begin{center}
\begin{tabular}{cc}
\includegraphics[width=6cm]{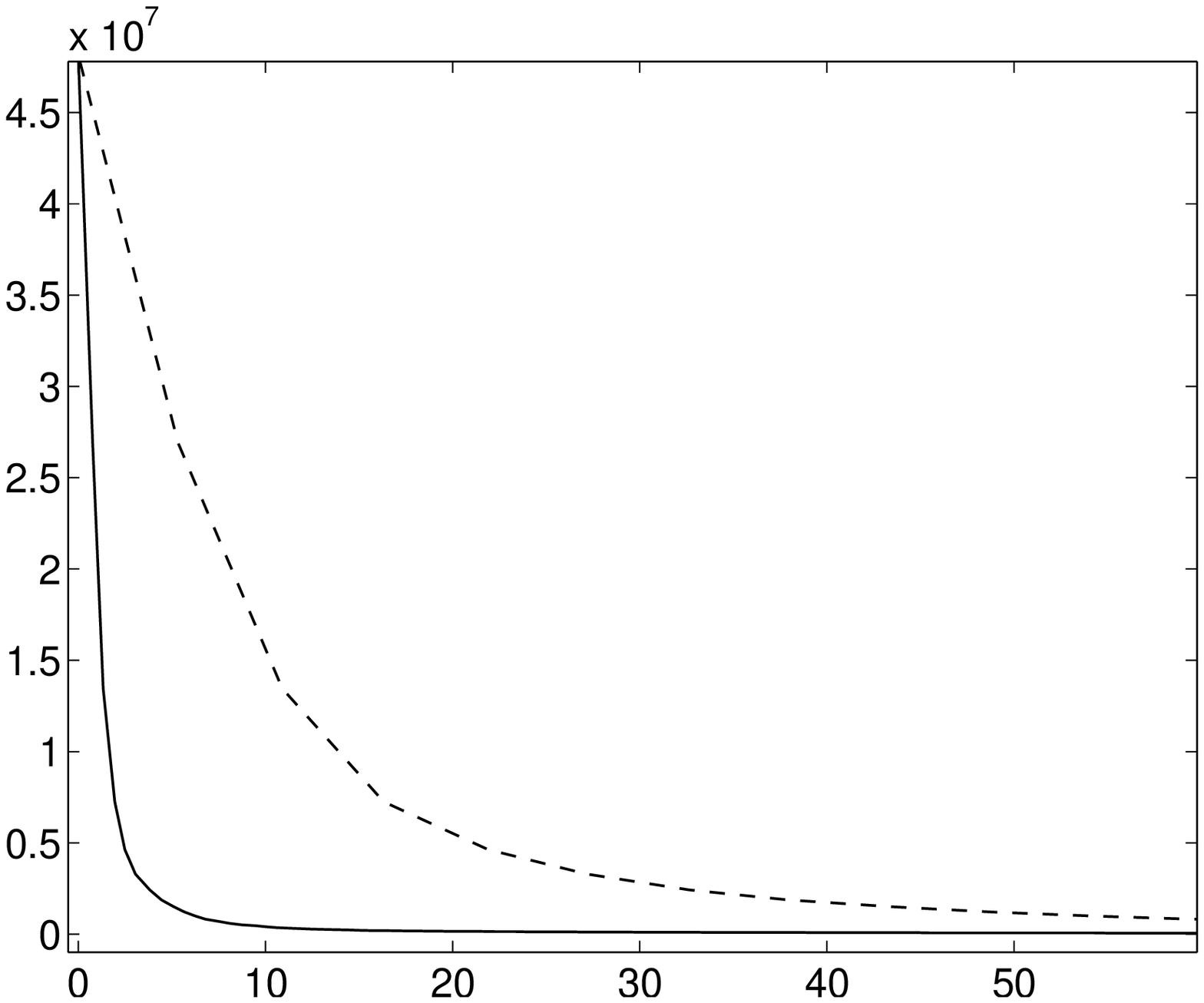}&\includegraphics[width=6cm]{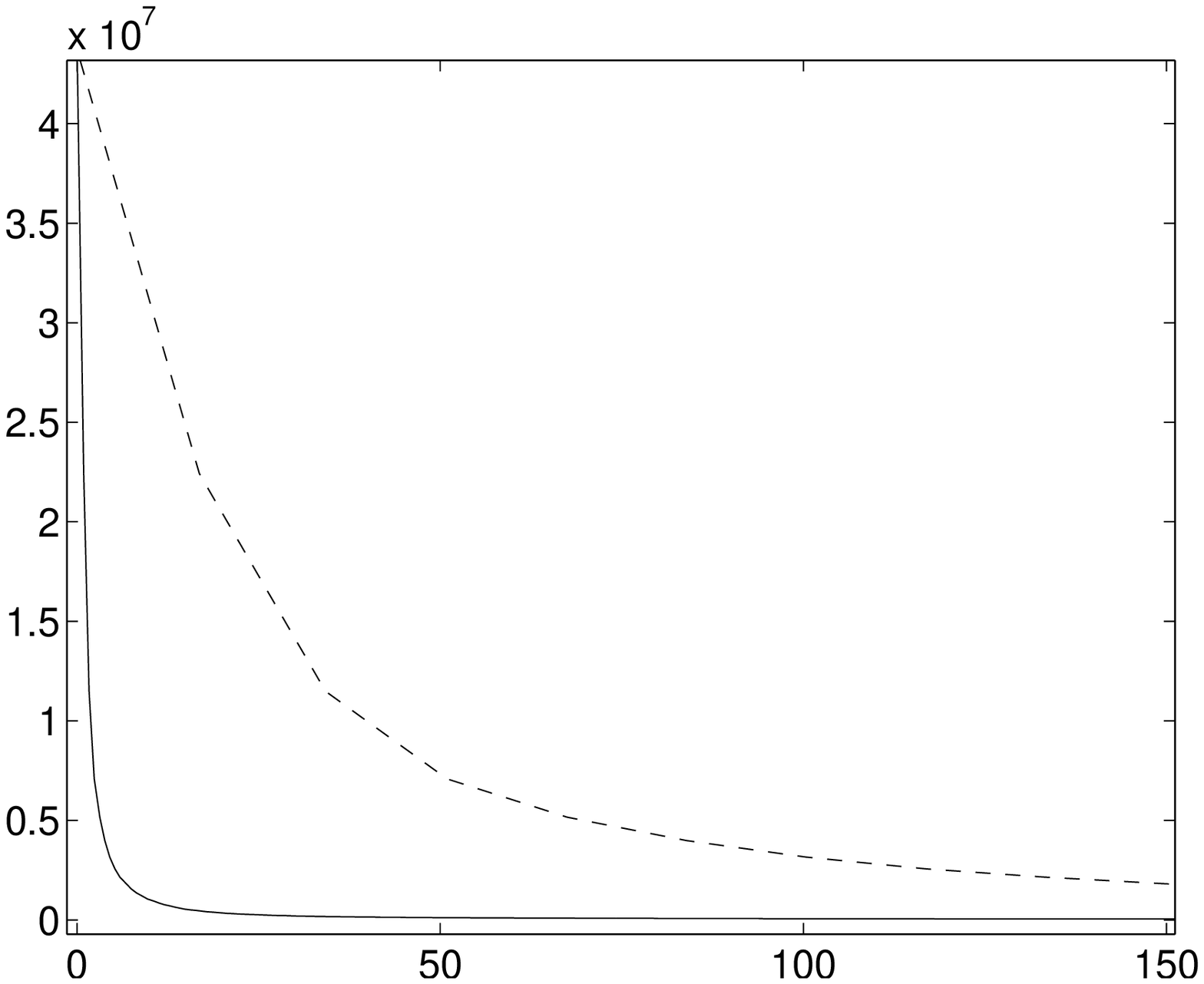} \\
\end{tabular}
\end{center}
\caption{Convergence profiles of the Algorithm~\ref{algo:PPA_frame} (dotted line) and Algorithm~\ref{algo:accPPA} (solid line) versus computation 
  time in seconds for a $3\times3$ uniform blur (left) and a $7\times7$ uniform blur (right) and a $128 \times128$-image. \label{fig:compCPU}}
\end{figure*}

Table \ref{tab:compVitessePPXA} gives the iteration numbers and the CPU times for the original PPXA algorithm and the proposed accelerated one
in order to reach convergence when considering different image sizes (``Sebal'': $N=128\times 128$, ``Peppers'': $N=256\times 256$ and ``Marseille'': $N=512\times 512$) and various kernel blur sizes.
The stopping criterion is based on the relative error between the objective function computed at the current iteration and at the previous one.\footnote{The relative error was evaluated based on Criterion \eqref{eq:minimization}
where the indicator function was discarded.}
 The stopping tolerance has been set to $10^{-3}$. These results have been obtained with an Intel Core2 6700, 2.66 GHz. 
The last line of Table~\ref{tab:compVitessePPXA} illustrates the gain in CPU-time when using Algorithm \ref{algo:accPPA}. 
Moreover, in Figure~\ref{fig:compCPU}, 
the mean square error on the image iterates
$\|F^{\top} (x_{n}- \hat{x})\|^2$ is plotted as a function of computation time, where $(x_n)_{n>0}$ denotes the sequence generated by Algorithm~\ref{algo:PPA_frame} or Algorithm~\ref{algo:accPPA}.

\begin{table*}[htbp]
\centering
\begin{tabular}{|c || c | c || c | c || c | c | }
\hline
Image size	& \multicolumn{2}{c||}{$128\times 128$}& \multicolumn{2}{c||}{$256\times 256$}& \multicolumn{2}{c|}{$512\times 512$}\\
\hline
(uniform) blur size& $3\times 3$&  $7\times 7$& $3\times 3$& $7\times 7$& $3\times 3$&  $7\times 7$\\
\hline
\hline
Iteration numbers & 30 & 50 & 41 & 50 & 50 & 50\\
\hline
CPU time (in second) & 117.2 & 633.0 & 411.7 & 1298 & 1458 & 4514 \\
\hline
CPU time - accelerated version (in second)& 13.53 & 29.82 & 60.59 & 89.48 & 263.6 & 405.0\\
\hline
Gain & 8.67 & 21.2 & 6.79 & 14.5 & 5.53 & 11.1 \\
\hline
\end{tabular}
\caption{Comparisons between PPXA and its accelerated version. \label{tab:compVitessePPXA}}
\end{table*}

It can be noticed that the larger the kernel blur size is, the higher the gain is. This is due to the fact that the number of proximity operators to compute increases with the kernel size.

\subsubsection{Comparison between different forms of total variation}
\begin{figure*}[htbp] 
\begin{center}
\begin{tabular}{c c}
\includegraphics[width=6cm]{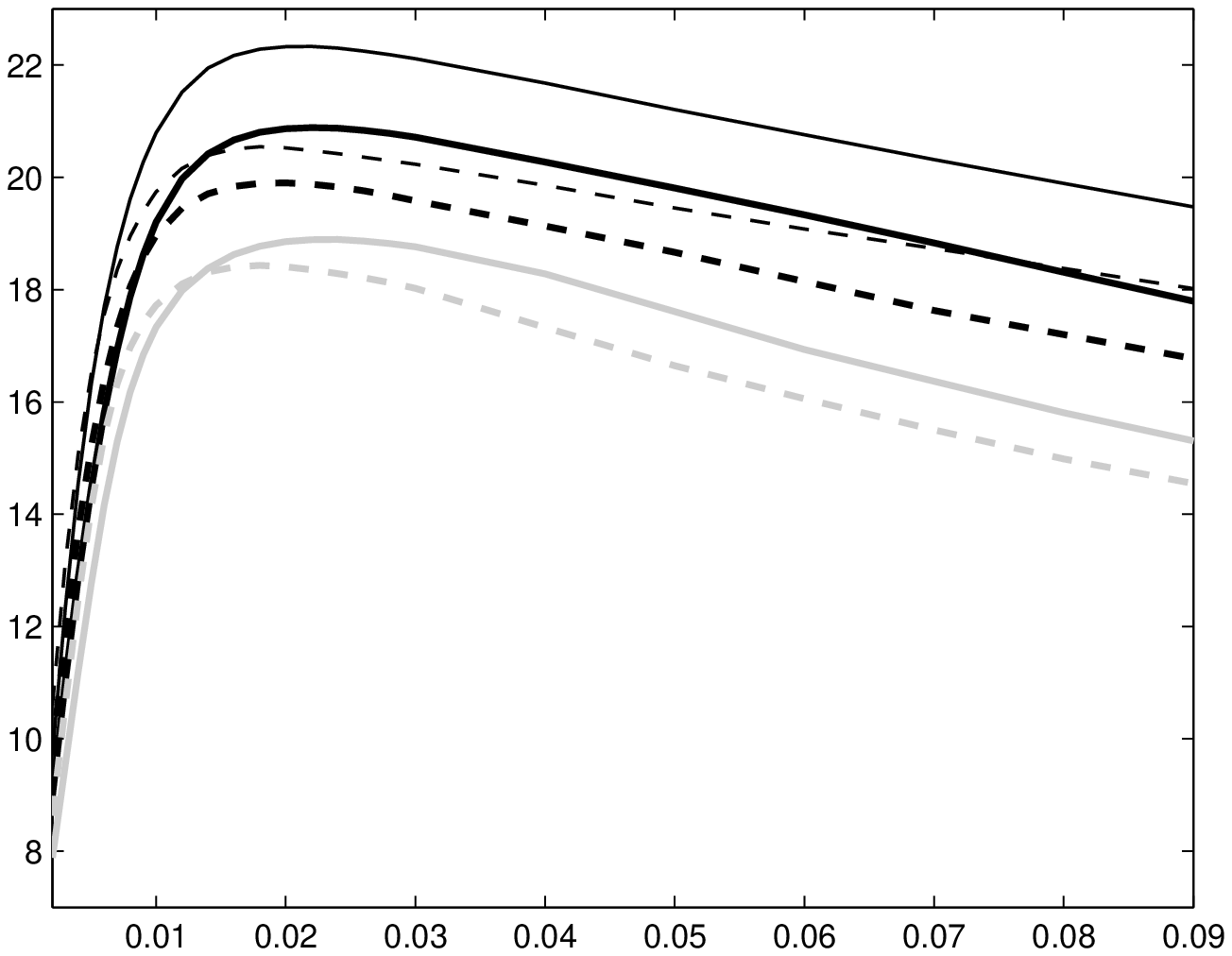} & \includegraphics[width=6cm]{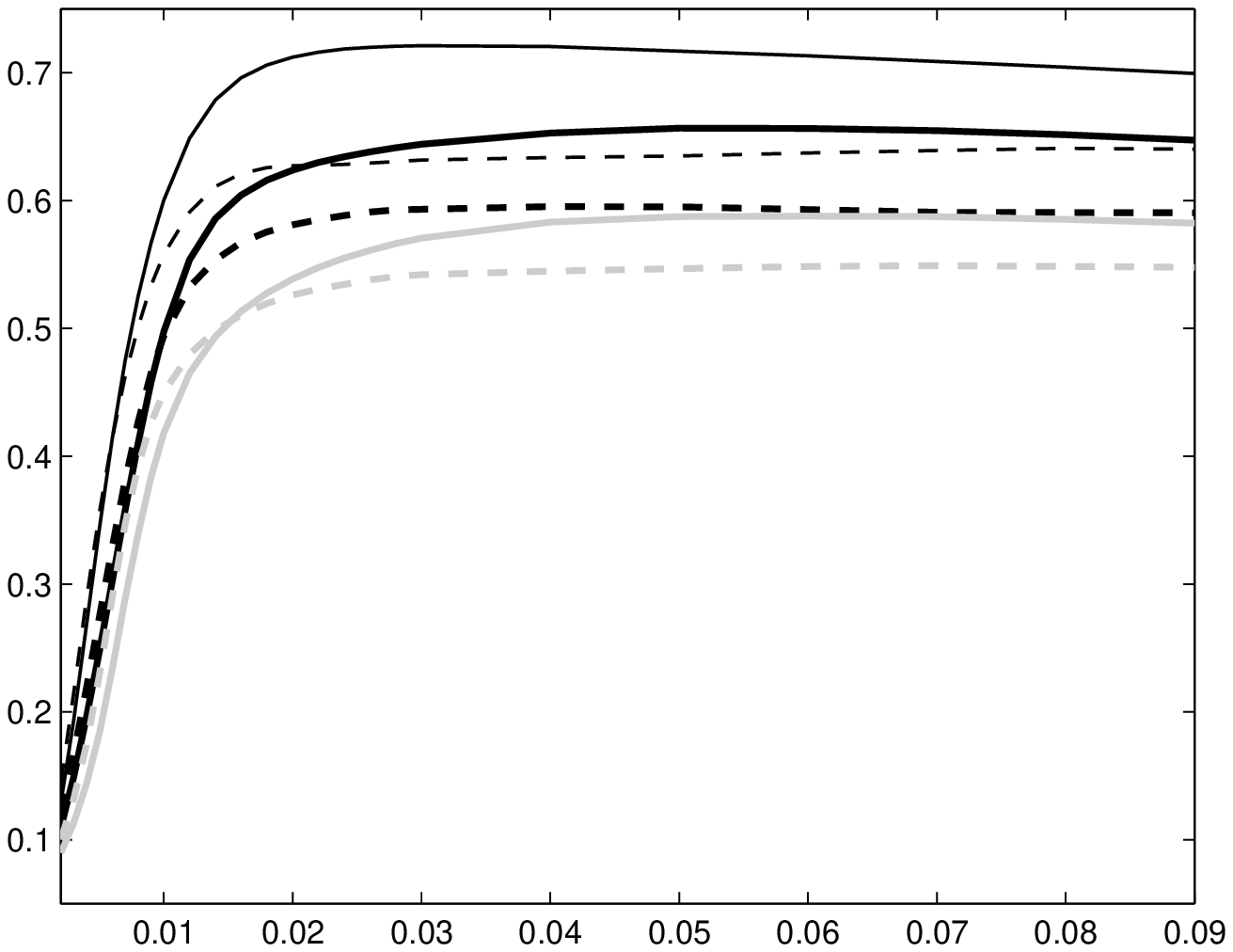}\\
\end{tabular}
\end{center}
\caption{SNR (left) and SSIM (right) for different total variation terms with respect to $\mu$. 
  Roberts: thin-black line, Finite difference: thick-black, Prewitt: thick-gray, $\rho_{\tv} = \sqrt{|\cdot|^2 + |\cdot|^2}$: solid line and $\rho_{\tv} = |\cdot| + |\cdot|$: dashed line \label{fig:compTV}}
\end{figure*}

In Section \ref{ss:proxtv}, we have introduced the proximity operator associated with discretized total variation functions and the possibility of choosing various filters has been mentioned. In Figure \ref{fig:compTV}, tests have been carried out on ``Peppers'' degraded by a $3\times 3$ uniform blur
and corrupted by Poisson noise with scaling parameter $\alpha = 0.1$.
We compare the restored images for different kinds of total
variation, in terms of Signal to Noise Ratio -- SNR 
and structural similarity measure -- SSIM 
  \cite{Wang_Z_2009_spm_MES_lioli}. The SSIM takes a value between -1
  to 1, the maximum value being obtained for two identical images. Each curve represents the resulting SNR and SSIM versus 
$\mu$ (the regularization parameter related to the total variation), for a given form of $\tv$ (i.e. a given filter associated with either an isotropic or anisotropic function $\rho_\tv$). A small wavelet regularization parameter $\vartheta=10^{-3}$ has been chosen in order to better illustrate the influence of the different $\tv$ forms on restoration quality.

It can be concluded from Figure \ref{fig:compTV} that the choice of the gradient filters and of the form (isotropic/aniso\-tropic) of $\rho_\tv$ has a significant influence on the restoration quality when the wavelet regularization is small. However, we also noticed in our numerical experiments that when the wavelet regularization parameter $\vartheta$ becomes larger, the choice of the $\tv$ form has a lower influence on the restoration quality provided that the regularization parameters are appropriately chosen.

\subsubsection{Boundary effects on restored images}

\begin{figure*}[htbp] 
\begin{center}
\begin{tabular}{c c}
\includegraphics[width=5cm]{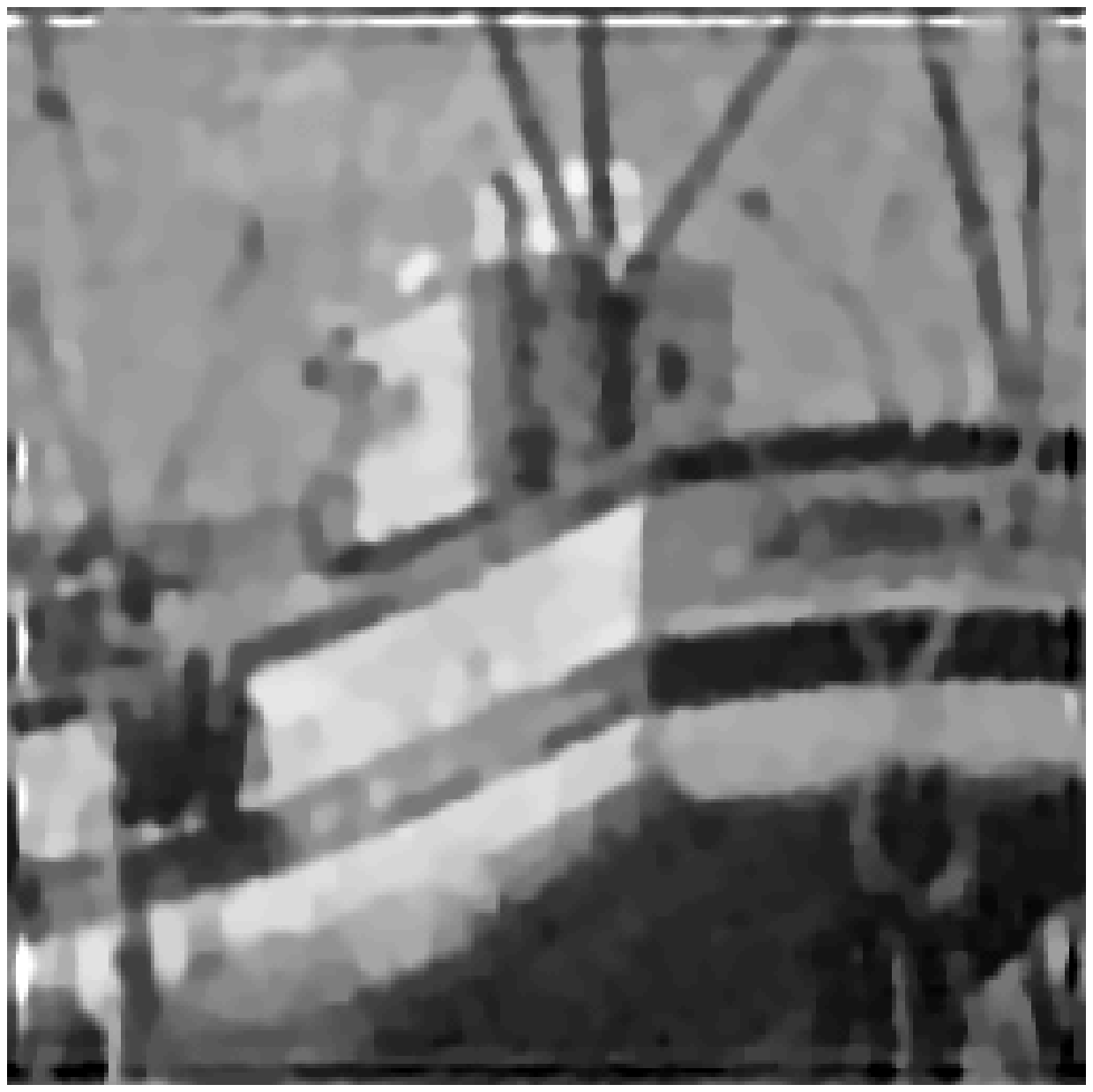} & \includegraphics[width=5cm]{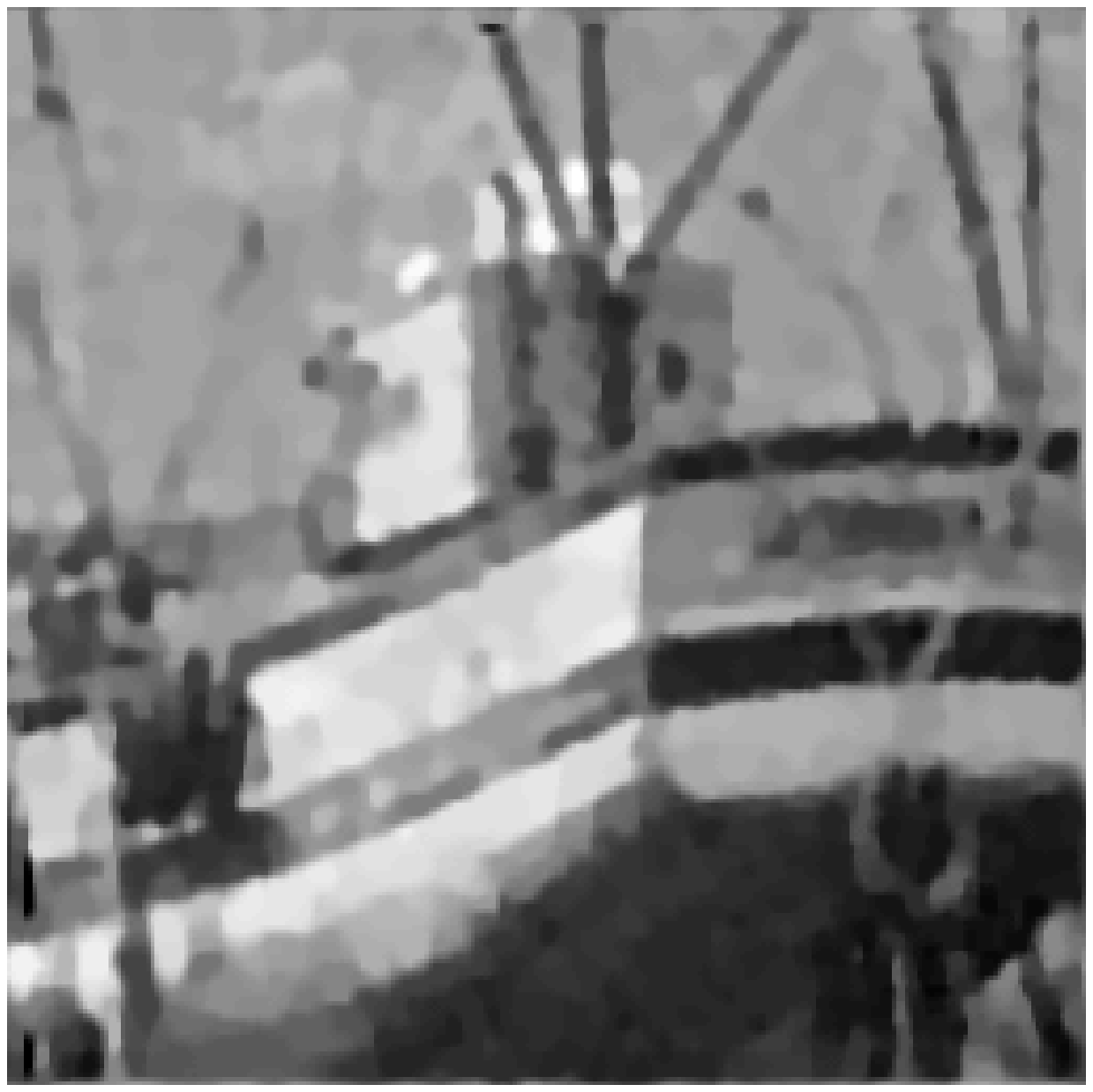}\\
\small{SNR = $16.9$ dB - SSIM = $0.62$} & \small{SNR = $17.7$ dB - SSIM = $0.64$}
\end{tabular}
\end{center}
\caption{Periodic (left) and zero-padded (right) restoration. \label{fig:comp_ZP_Per}}
\end{figure*}

This section illustrates the influence of boundary effect processing. More precisely, we degraded an extended version of ``Boat'' image by a $7\times 7$ uniform blur, and the resulting blurred image was cropped to create an image of size $256\times 256$.  As a consequence, the boundary values are functions of pixel locations which are no longer present in the blurred image. The scaling parameter associated with Poisson noise is $\alpha=0.5$. The objective is then to restore the image (which was centered) by using one of the convolution
models discussed in Section~\ref{ss:proxconv}, namely either
a periodic convolution or a convolution with zero-padding.
Visual and quantitative results are given in Figure~\ref{fig:comp_ZP_Per}.

As it can be noticed from this figure, the periodic convolution model introduces significant boundary artefacts unlike the convolution with
zero-padding. The results obtained when considering ``Peppers'' led to the same conclusion. For ``Sebal'', zero-padding or periodic models provided similar results.

\subsubsection{Decimated convolution}
\begin{figure*}[htbp]
\centering
\begin{tabular}{c c c}
\includegraphics[width=5cm]{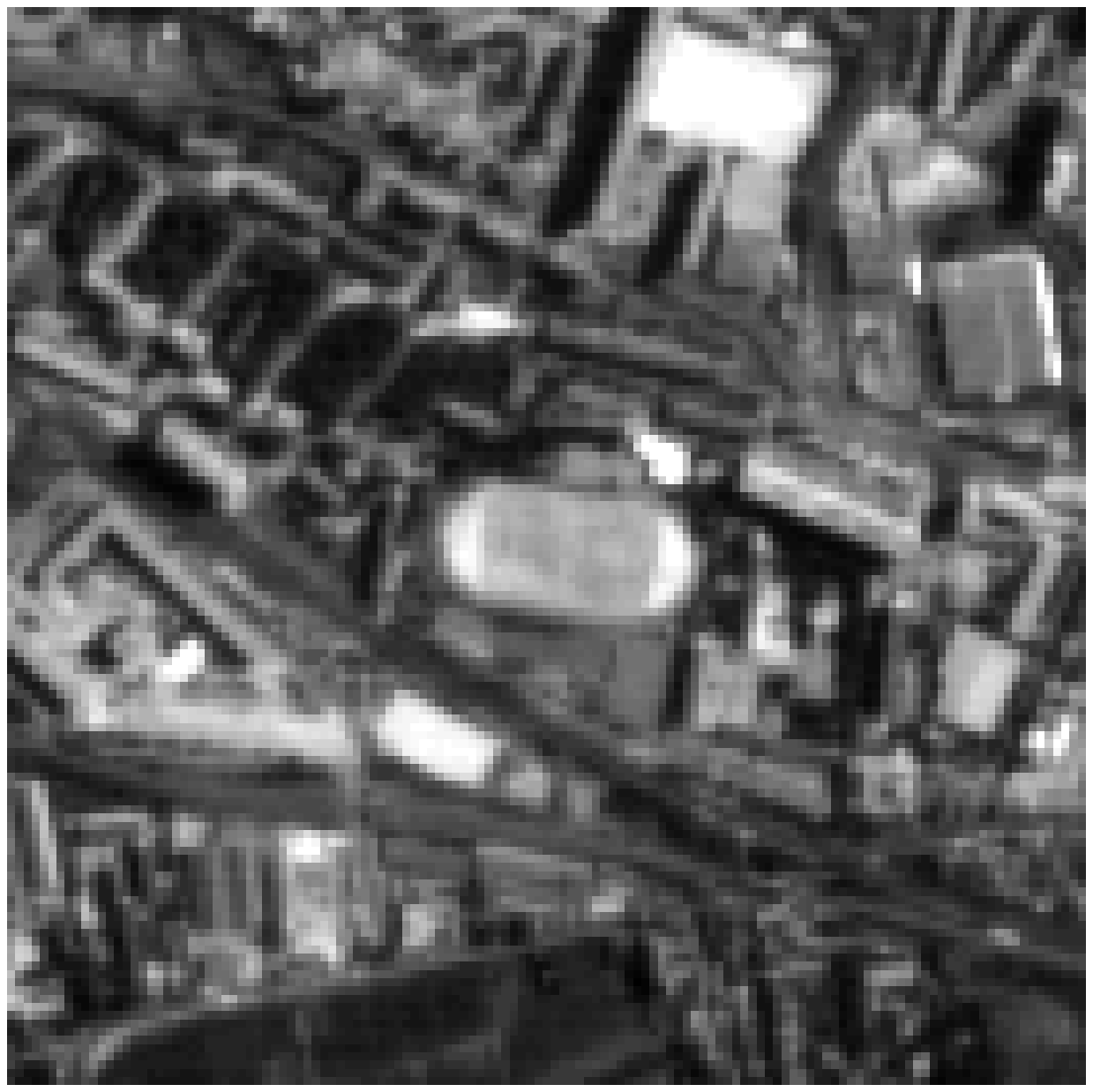}& \includegraphics[width=5cm]{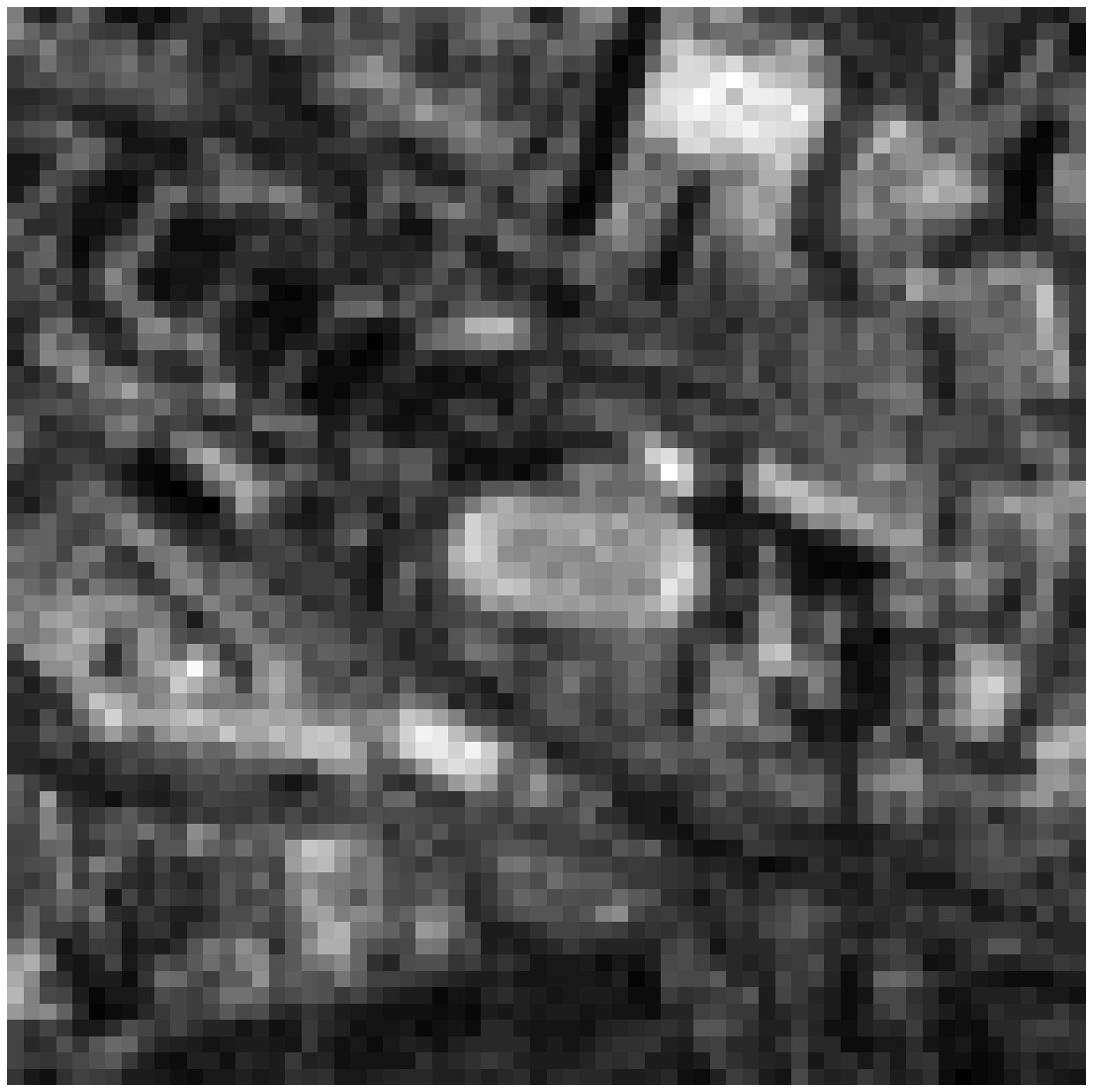} & \includegraphics[width=5cm]{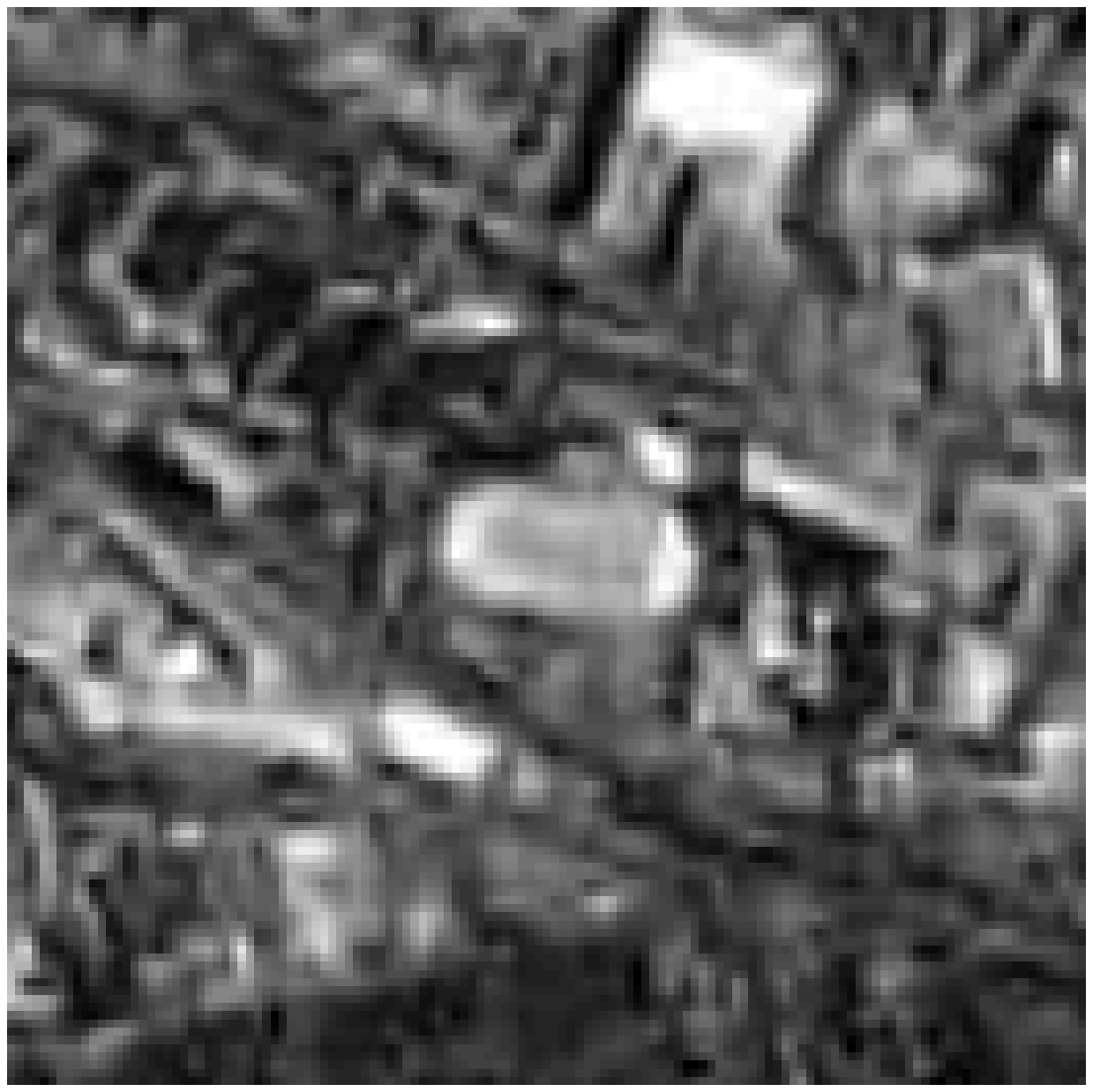}\\
\small{Original}&\small{Degraded}&  \small{Restored $(\vartheta=1 ,\mu=10^{-3})$} \\ 
 & &  \small{SNR = 16.1 - SSIM = 0.79} \\ 

\end{tabular}
\caption{Restoration results for ``Spot'' image.\label{fig:spot}}
\end{figure*}

We now present experimental results for a  $256\times 256$ SPOT image degraded by a uniform decimated blur with a kernel size $Q = 3\times 3$ and a decimated factor $d=2$. The scaling parameter of the Poisson noise is equal to $\alpha=1$. Due to the structure of the degradation operator, the data fidelity is splitted in a sum of $I=4$ functions. The results are presented in Figure~\ref{fig:spot} where the good behaviour of the model can be observed.

\subsubsection{Influence of each regularization term}
We now present numerical and visual results for the different kinds of regularization when a generalized Kullback-Leibler divergence is used as a data fidelity term. 
This experiment allows us to compare the hybrid
  regularization with existing approaches based on a wavelet-frame
  \cite{Figueiredo_M_2010_t-ip_restoration_piado} regularization or
  a total variation regularization
  \cite{Setzer_S_2009_Deblurring_pibsbt,Figueiredo_M_2010_t-ip_restoration_piado}. The latter
  regularized solutions can be computed either by using augmented
  Lagrangian techniques
  \cite{Setzer_S_2009_Deblurring_pibsbt,Figueiredo_M_2010_t-ip_restoration_piado}
  or with our splitting approach (by setting $\vartheta = 0$ or $\mu =
  0$). In our experiments, the computation time of the two approaches was observed
  to be similar. 
Note that comparisons performed in
\cite{Figueiredo_M_2010_t-ip_restoration_piado} led to the conclusion
that the wavelet-frame regularization is quite competitive with respect to other
existing restoration methods
\cite{Dupe_FX_2008_ip_proximal_ifdpniusr}.

\begin{figure*}[htbp]
\centering
\begin{tabular}{c c}
\includegraphics[width=6cm]{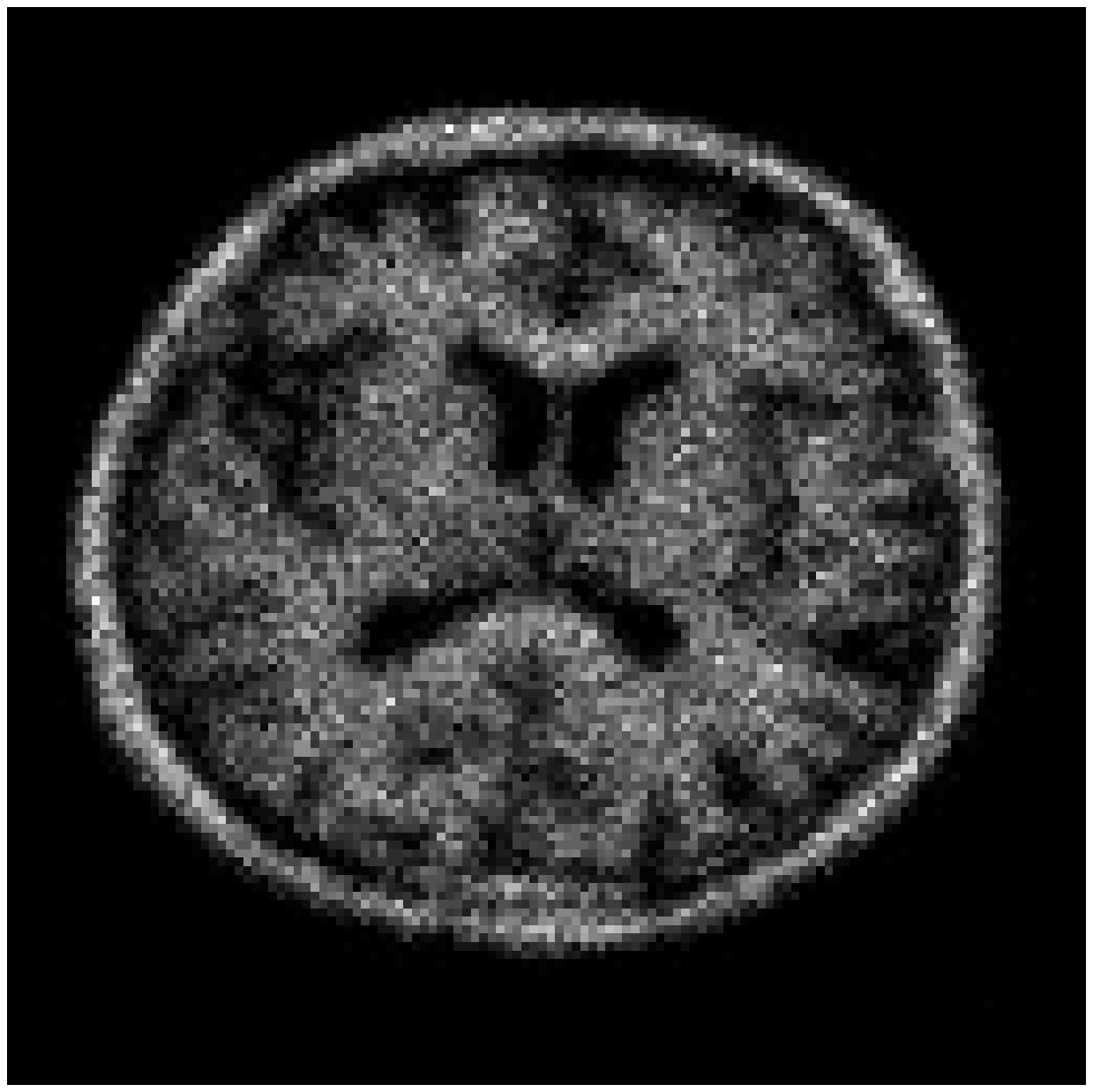}& \includegraphics[width=6cm]{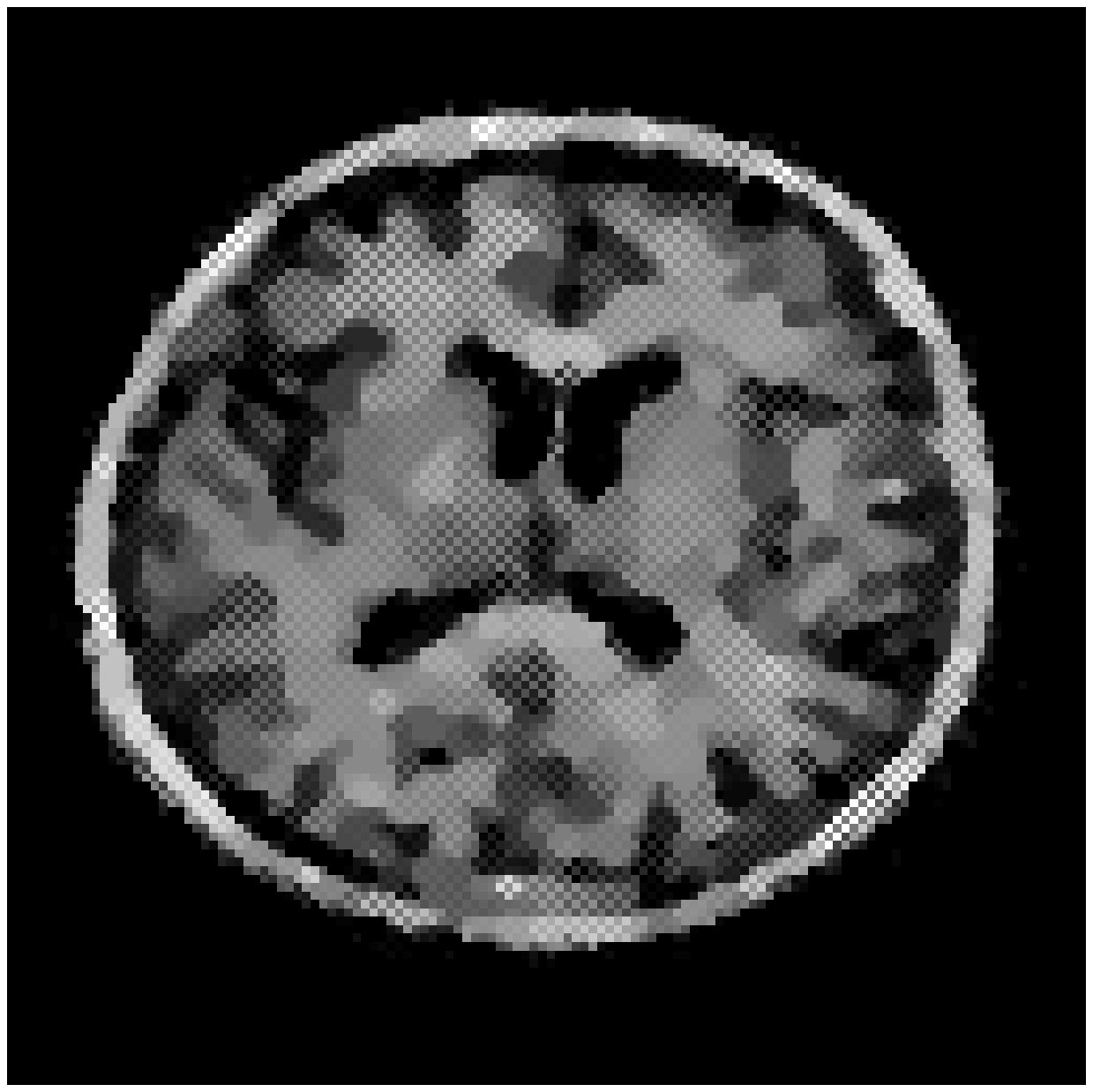} \\
\small{Degraded, $\alpha = 0.1$ and uniform blur $3 \times 3$}&  Total variation regularization \\
\small{SNR = 8.88 dB - SSIM = 0.69 }
&  \small{SNR = 11.2 dB - SSIM = 0.79}\\ 
\includegraphics[width=6cm]{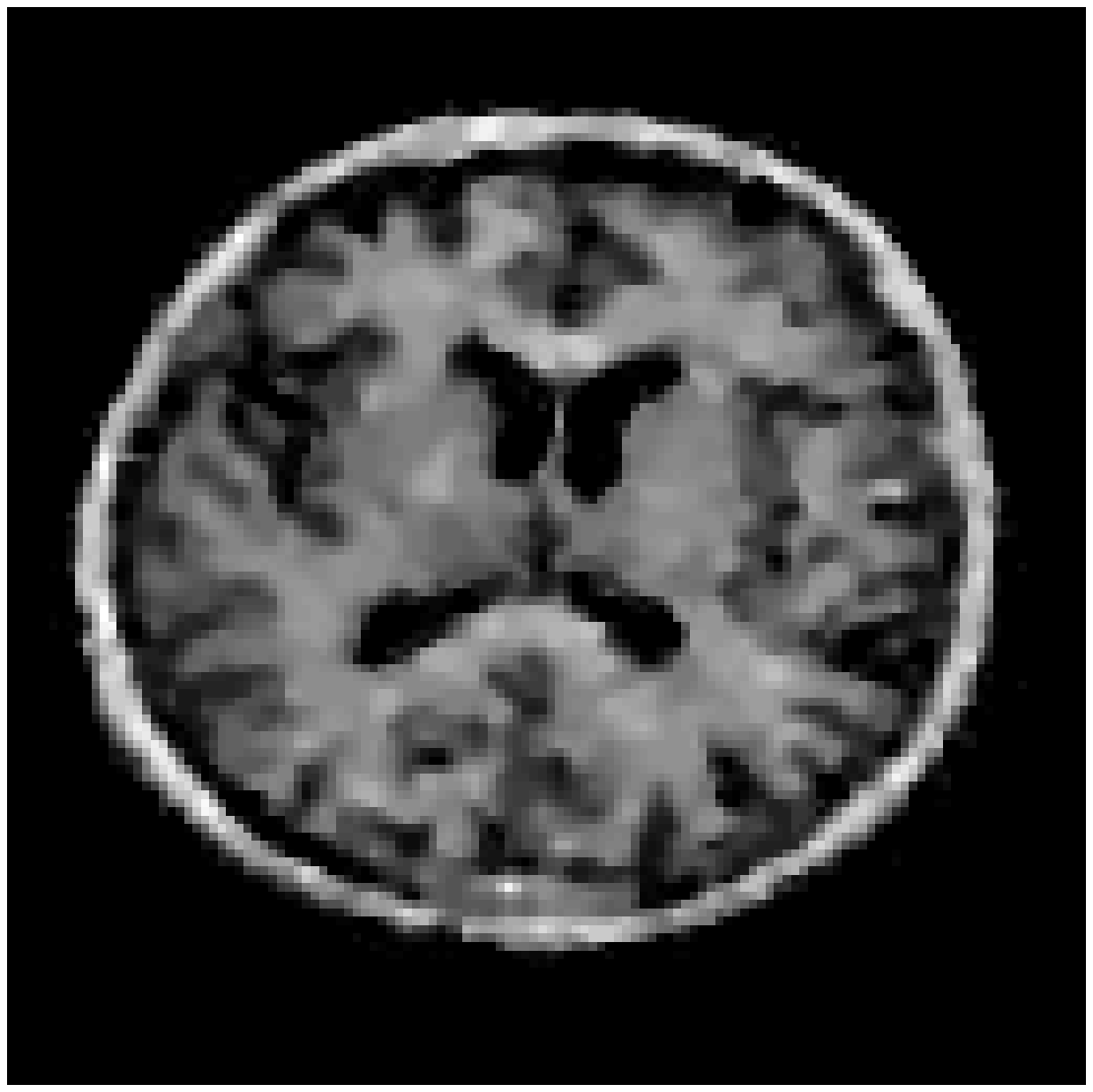}&  \includegraphics[width=6cm]{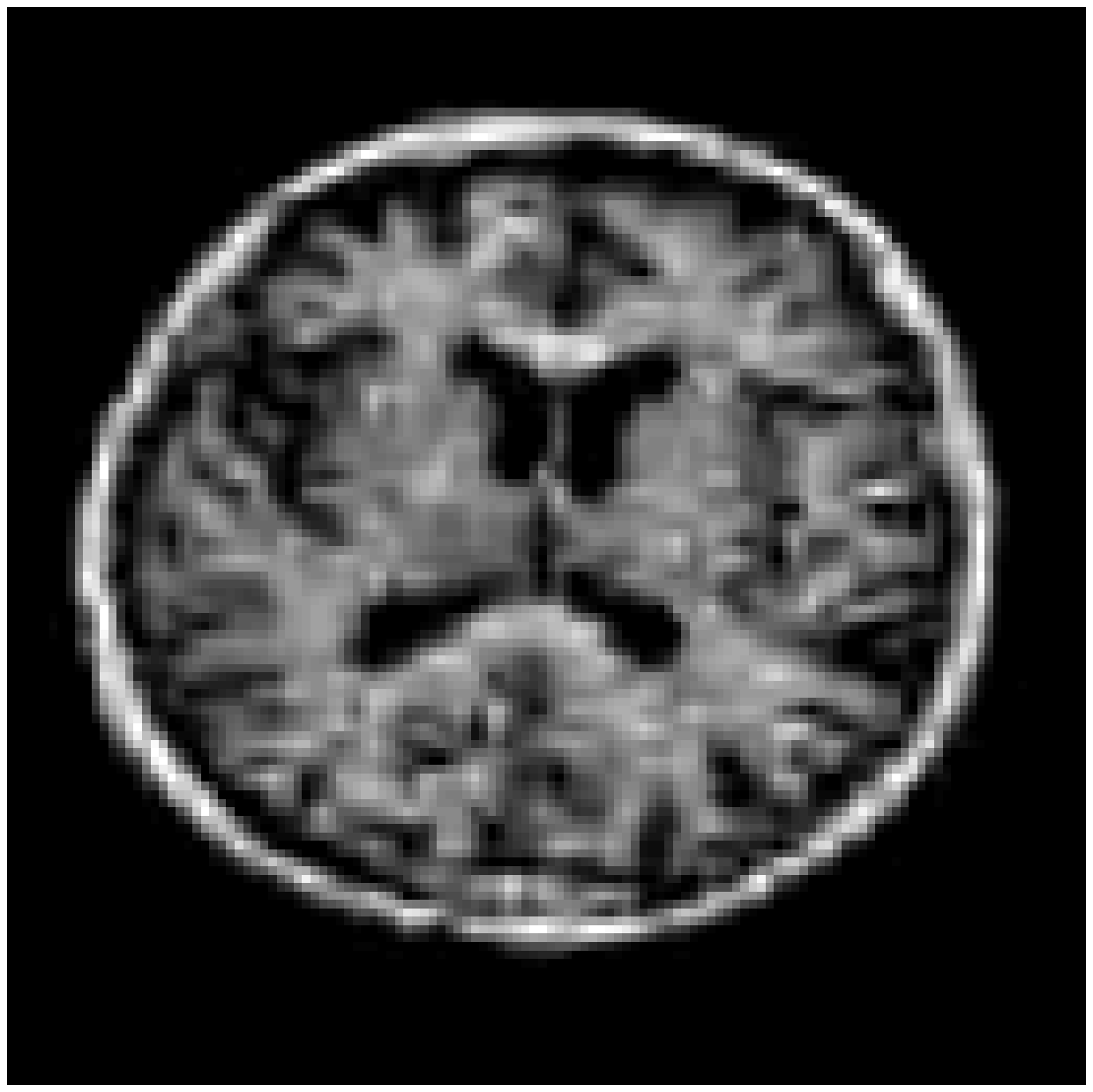} \\
Hybrid regularization  \small{$(\vartheta=0.09 ,\mu=0.006)$} & Wavelet-frame regularization \\
\small{SNR= 12.4 dB - SSIM= 0.85} & \small{SNR= 11.7 dB - SSIM= 0.83}\\
\end{tabular}
\caption{Restoration results for ``Sebal'' image.\label{fig:sebal}}
\end{figure*}

\begin{figure*}[htbp]
\centering
\begin{tabular}{c c}
\includegraphics[width=6cm]{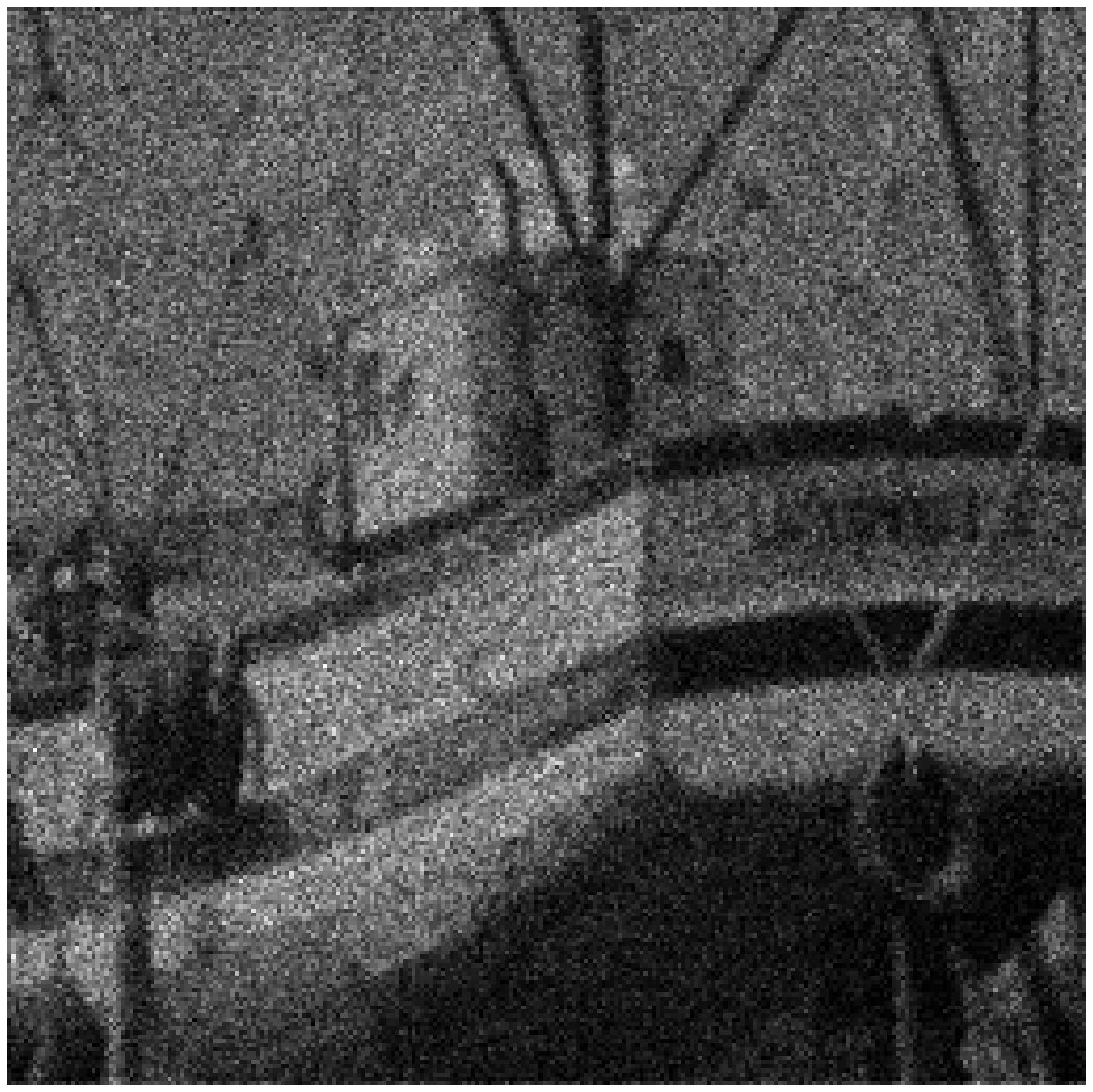}& \includegraphics[width=6cm]{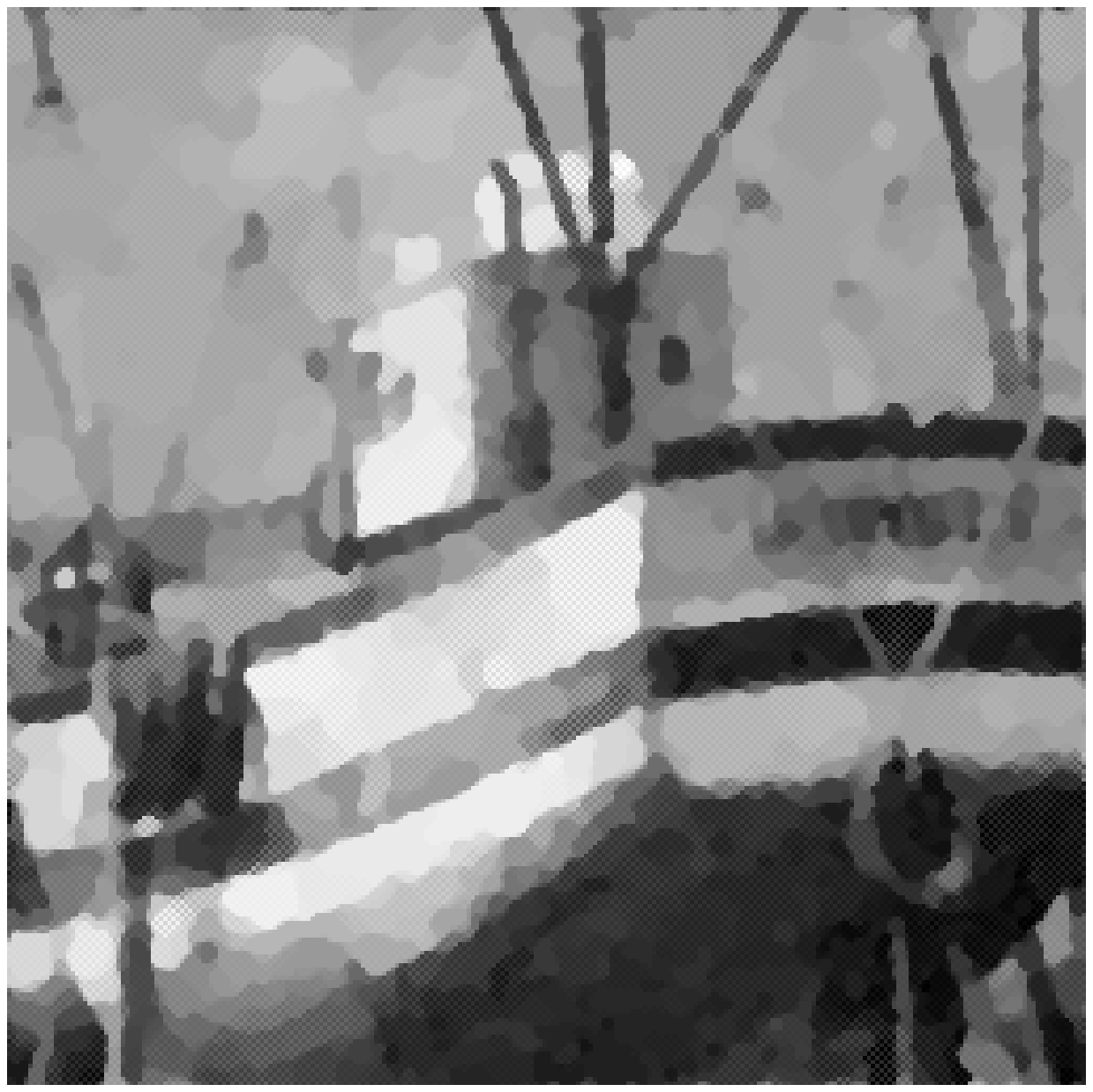} \\
\small{Degraded, $\alpha = 0.1$ and uniform blur $3 \times 3$}&  Total variation regularization \\
\small{SNR = 11.2 dB - SSIM = 0.27 }&  \small{SNR = 17.8 dB - SSIM = 0.60 } \\
\includegraphics[width=6cm]{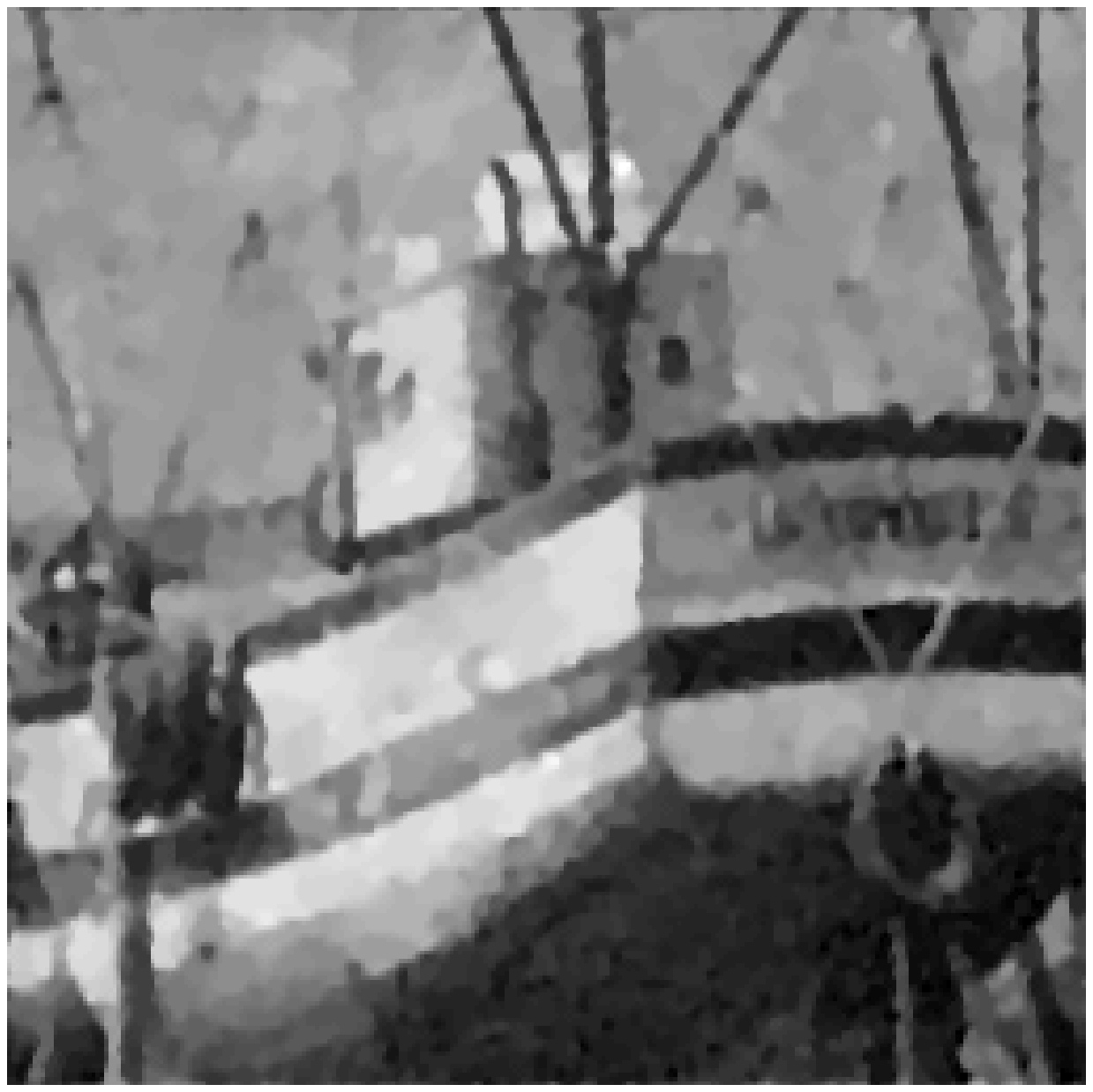}&  \includegraphics[width=6cm]{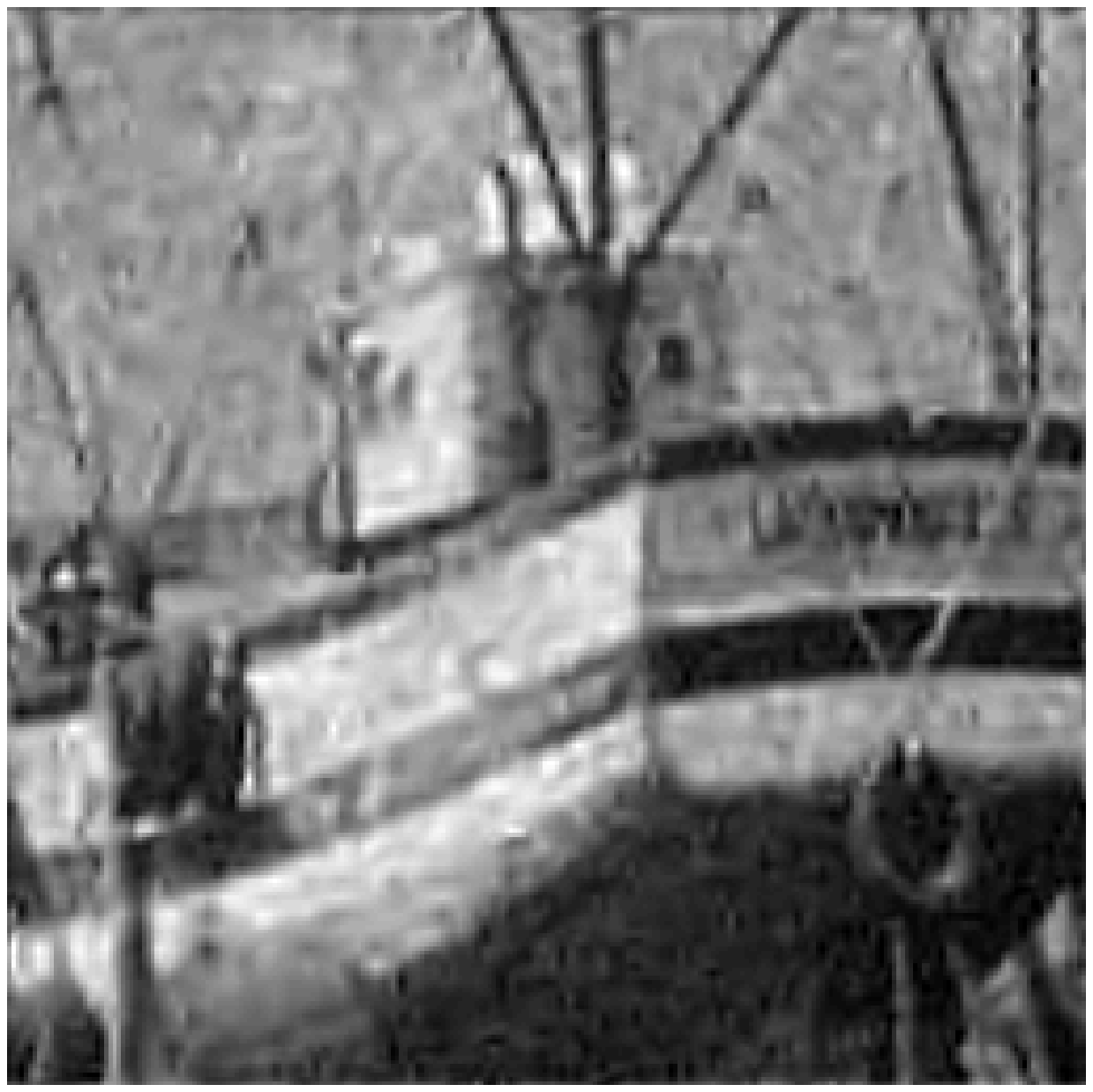} \\
Hybrid regularization  \small{$(\vartheta=0.06 ,\mu=0.011)$} 
&  Wavelet-frame regularization \\
\small{SNR=18.8 dB - SSIM=0.67} & \small{SNR=18.0 dB - SSIM= 0.62}\\
\end{tabular}
\caption{Restoration results for ``Boat'' image.\label{fig:boat}}
\end{figure*}

\begin{figure*}[htbp]
\centering
\begin{tabular}{c c}
\includegraphics[width=6cm]{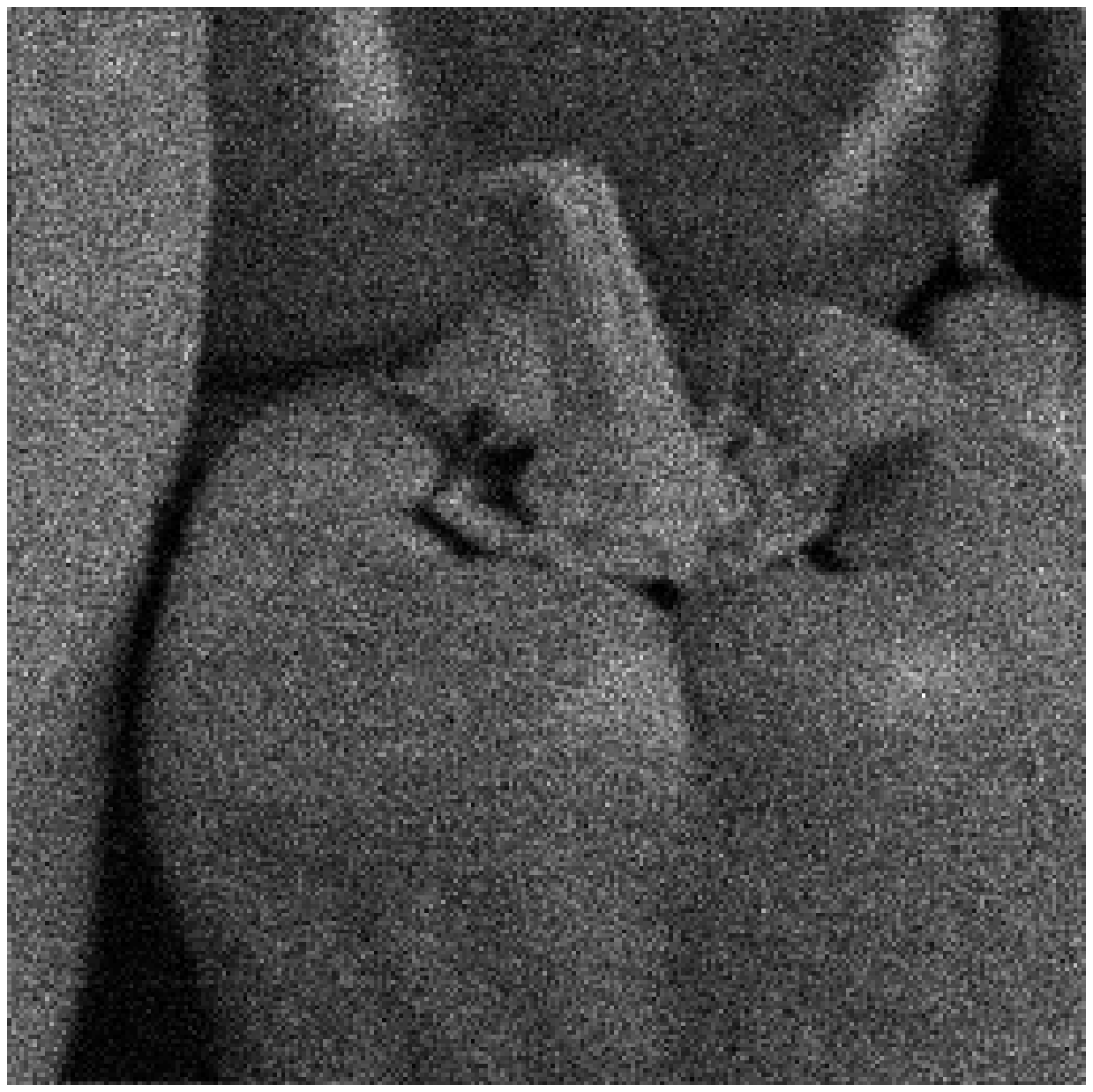} &  \includegraphics[width=6cm]{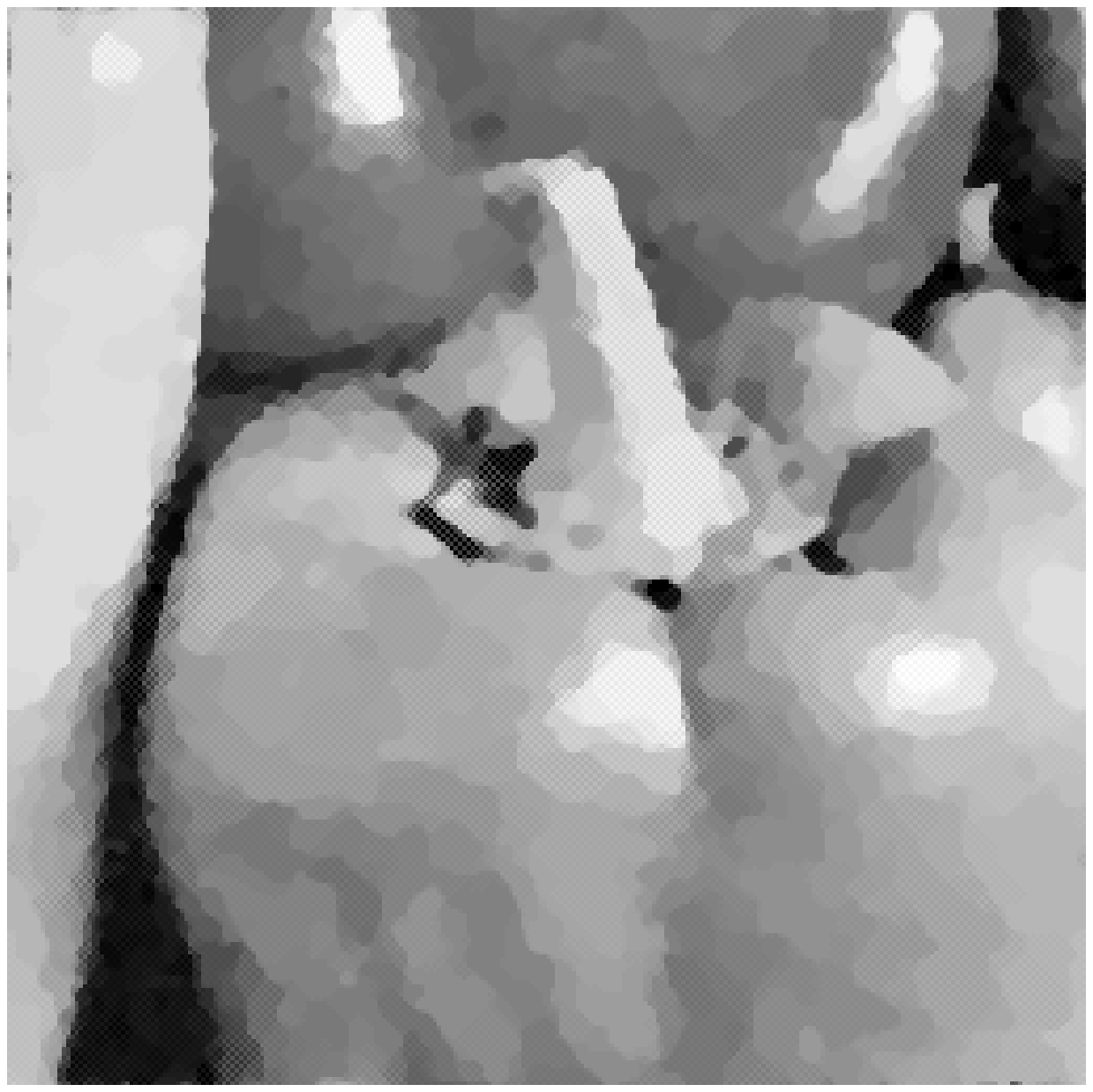}\\
\small{Degraded, $\alpha = 0.1$ and uniform blur $3 \times 3$}&  Total variation regularization \\
\small{SNR = 11.4 dB - SSIM = 0.16}&  \small{SNR = 22.1 dB - SSIM = 0.69}\\
\includegraphics[width=6cm]{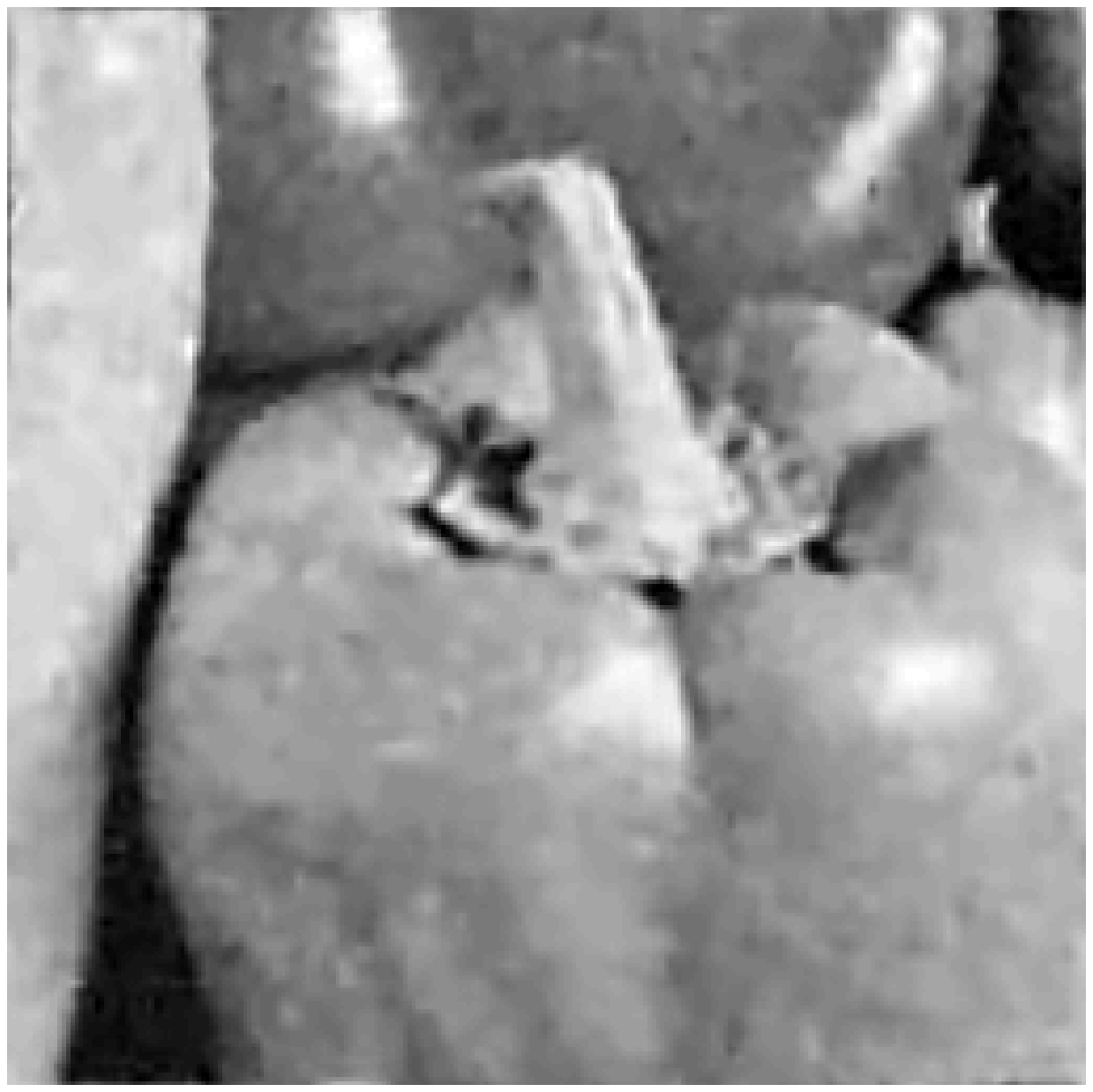} & \includegraphics[width=6cm]{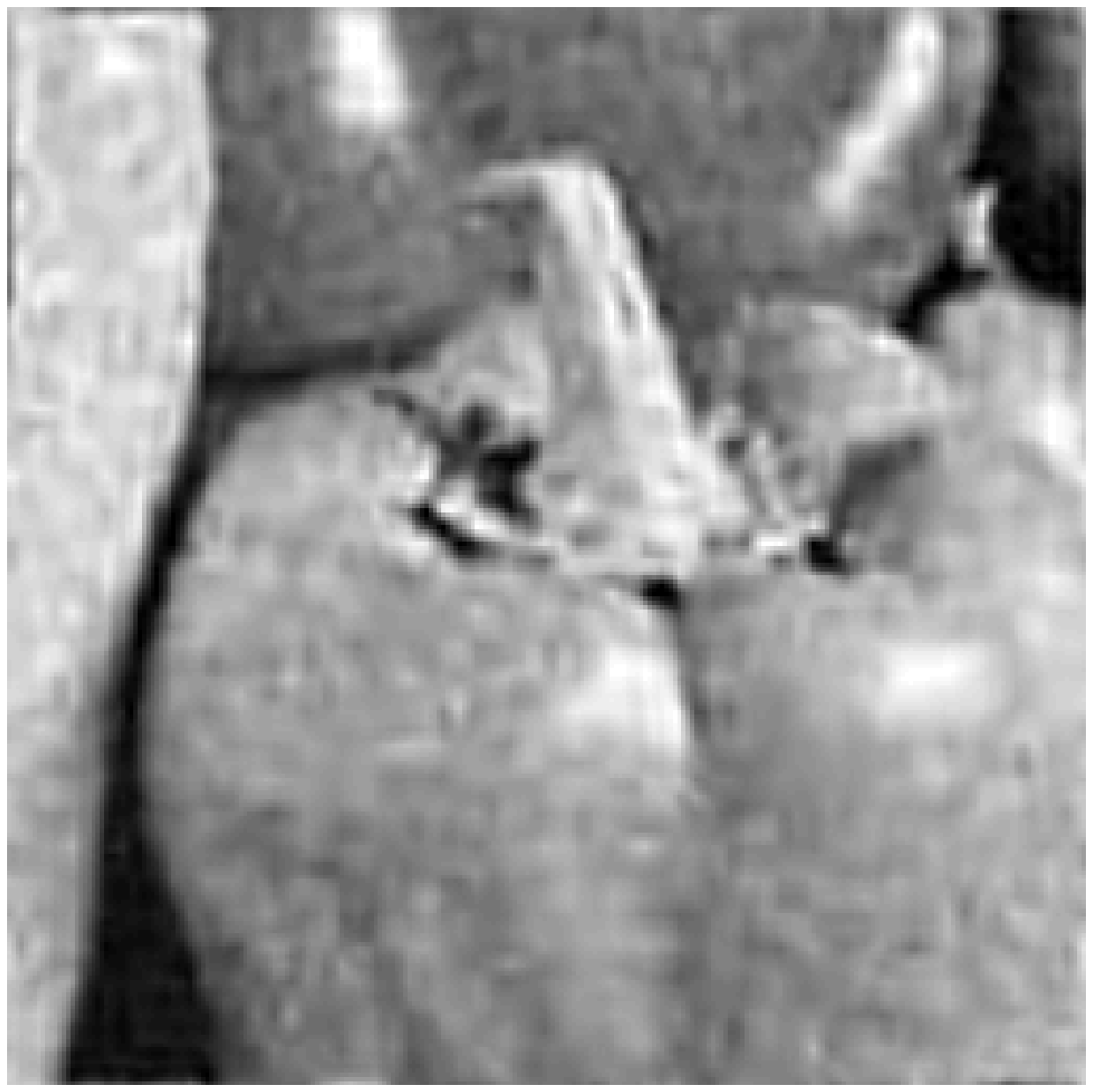}\\
Hybrid regularization \small{$(\vartheta=0.2 ,\mu=0.006)$} 
&   Wavelet-frame regularization\\ 
\small{SNR = 22.7 dB - SSIM = 0.74} & \small{SNR= 21.4 dB - SSIM = 0.68 }\\
 \end{tabular}
\caption{Restoration results for ``Peppers'' image.\label{fig:peppers}}
\end{figure*}

In the images displayed in Figures \ref{fig:sebal}, \ref{fig:boat}, and
\ref{fig:peppers}, one can observe the
artefacts related to the wavelet regularization, the staircase effects
which are typical of the total variation penalization, some checkerboard patterns resulting of the chosen gradient discretization, and also the
benefits which can be drawn from the use of a hybrid regularization. 

Similarly to \cite{Figueiredo_M_2010_t-ip_restoration_piado},
  the values of $\mu$ and $\vartheta$ have been adjusted  so as to maximize the SNR. Optimizing the hyperparameters manually as we did
is a common practice in imaging applications, especially when a data
set of test images having similar characteristics as the one to be
restored (medical images, satellite images,...)
is available. Automatic methods for the optimization of the hyperparameters can also be found in the literature such as cross-validation \cite{Galatsanos_N_1992_tip_methods_crpenvirtr}, Stochastic EM \cite{Delyon_B_1999_annals-statistics_convergence_savema},
MCMC
\cite{Robert_C_2004_monte_csm,Chaari_L_2010_j-tsp_hierarchical_bmfr}
or Stein-based methods
\cite{Ramani_S_2008_tip_monte_csbborpgda}. These automatic procedures
often are relatively intensive. They will not be addressed in this paper due to the lack of space.

\section{Conclusion}\label{se:conclu}
A new convex regularization approach to restore data degraded by a
(possibly decimated) convolution operator and a non necessarily additive noise has been proposed. The main advantages of the method are (i) to deal directly with the ``true'' noise likelihood (i.e. the Kullback-Leibler divergence in the case of Poisson noise) without requiring any approximation of it; (ii) to permit the use of sophisticated regularization functions, e.g. one promoting sparsity in a wavelet frame domain and a total variation penalization. In addition, 
the proposed algorithm has a parallel structure which makes it easily implementable on multicore architectures. 
Numerical and visual results demonstrate the effectiveness of the proposed approach. One can note that, even if the paper is devoted to the case of convolutive operators, this approach could be generalized to more general  linear operators.

Note that the primal-dual approaches
  \cite{Zhu_M_2007_report_2008_efficient_pdhgatvir,Zhang_X_2010_jsc_unified_pdafbbi,Esser_E_2010_report_general_fcfopdatvm,Stezer_S_2010_cms_infimal_crdl1tf,Chambolle_A_2010_first_opdacpai}
  can offer alternative solutions to the ones developed in this
  paper. However, one of the advantages of PPXA is that it easily leads to
efficient parallel implementations


\section{Proof of Proposition \ref{p:proxcomp}} \label{ap:proxcomp}

Since $D=L L^{\top}$ is the matrix associated with  a bijective operator, $L$ is associated with a surjective one
and $\dom\Upsilon \neq \emp \Rightarrow \dom( \Upsilon \circ L)\neq \emp$.
This allows us to conclude that $h = \Upsilon \circ L$ is a function
of $\Gamma_0(\RR^X)$.\\
To calculate the proximity operator of $h$, we now come back to the definition of this operator. We have thus, for every $w\in \RR^X$,
\begin{equation}
\prox_h w = \arg \min_{v \in \RR^X} \frac{1}{2} \Vert v-w\Vert^2 + \Upsilon(Lv).
\end{equation}
We can write any vector $v\in \RR^X$ as a sum of an element $L^{\top}t \in \ran L^{\top}$ and $v_{\perp}\in (\ran L^{\top})^{\perp}= \ker L$. 
We have then 
$Lv = LL^{\top}t  = Dt$. Similarly, we can write $w = L^{\top}u+w_{\perp}$
where $u\in \RR^Y$ and $w_{\perp}\in\ker L$.
So, $\prox_h w$ can be determined by finding
\begin{multline}
\min_{(t,v_{\perp}) \in \RR^Y\times \RR^X} \frac{1}{2} \Vert L^{\top}t + v_{\perp} - L^{\top}u-w_{\perp}\Vert^2 + \Upsilon(Dt)\\
= \min_{(t,v_{\perp}) \in \RR^Y\times\RR^X} \frac{1}{2} \Vert L^{\top}(t-u)\Vert^2 + \frac{1}{2}\Vert v_{\perp} -w_{\perp}\Vert^2 +  \Upsilon(Dt).
\end{multline}
This yields 
\begin{equation}
v_{\perp}=w_{\perp}= w-L^{\top}u
\label{eq:yperp}
\end{equation}
and it remains to find
\begin{equation}
\min_{t\in \RR^Y} \frac{1}{2} \Vert L^{\top}(t-u)\Vert^2  +  \Upsilon(Dt) =\min_{t\in \RR^Y} \frac{1}{2} \scal{D(t-u)}{(t-u)}+  \Upsilon(Dt).
\end{equation}
By using the separability of $\Upsilon$, this is equivalent to finding
\begin{equation}
\min_{t\in\RR^Y}  \sum_{m=1}^Y \frac{1}{2} \Delta_m(\scal{t}{o_m}  -
\scal{u}{o_m} )^2 +  \psi_m(\Delta_m \scal{t}{o_m} ).
\end{equation}
It can be deduced from \cite[Lemma~2.6]{Combettes_PL_2005_mms_Signal_rbpfbs} that, for every $m \in \{1, \ldots, Y \}$,
\begin{align}
\scal{t}{o_m} &= \prox_{\frac{1}{\Delta_m}\psi_m(\Delta_m \cdot)}(\scal{u}{o_m})\nonumber\\
&= \frac{1}{\Delta_m}\prox_{\Delta_m\psi_m}(\Delta_m\,\scal{u}{o_m}),
\end{align}
which, according to \cite[Proposition~2.10]{Chaux_C_2007_ip_variational_ffbip}, leads to
\begin{align}
\label{e:proxz}
t &= D^{-1}\sum_{m=1}^Y
\prox_{\Delta_m\psi_m}(\Delta_m\scal{u}{o_m})\,
o_m \nonumber \\
&= 
D^{-1}\prox_{D\Upsilon}(Du).
\end{align}
Altogether, \eqref{eq:yperp} and \eqref{e:proxz}
yield  
\begin{equation*}
v = w + L^{\top} \big(D^{-1}\prox_{D\Upsilon}(Du) -u\big).
\end{equation*}
In addition, since $L^{\top}u$ is the projection of $w$ onto $\ran L^{\top}$,
$u = (LL^{\top})^{-1} Lw = D^{-1} Lw$ and \eqref{e:pfL} follows.\\

\section{Proof of Proposition \ref{p:proxtv}} \label{ap:proxtv}

By using the proximity operator definition \eqref{eq:prox}, 
\begin{equation*}
\pi =  \prox_{\mu \tv_{p_1,p_2}}(y) 
\end{equation*}
minimizes
\begin{align}
\frac{1}{2}\Vert \pi - &y\Vert^2 + \mu \tv_{p_1,p_2}(\pi) \nonumber\\
 =&  \;  \frac{1}{2} \sum_{(n_1,n_2)\in \mathbb{B}} 
(\pi_{n_1,n_2}-y_{n_1,n_2})^2\nonumber\\
&+ \sum_{n_1 = 0}^{\lfloor\frac{N_1-p_1}{P_1}\rfloor-1}\sum_{n_2 = 0}^{\lfloor\frac{N_2-p_2}{P_2}\rfloor-1}\Big\{\frac{1}{2}\Vert \Pi_{n_1,n_2}^{p_1,p_2} - Y_{n_1,n_2}^{p_1,p_2}\Vert^2_{\rm F}+ \mu\; \rho_{\rm tv}\big( \trace(\hh^\top \Pi_{n_1,n_2}^{p_1,p_2}),\trace(\vv^\top \Pi_{n_1,n_2}^{p_1,p_2})\big)\Big\}
\end{align}
where $\|\cdot\|_{\rm F}$ denotes the Frobenius norm, for every $(n_1,n_2) \in \{0,\ldots,\lfloor\frac{N_1-p_1}{P_1}\rfloor-1\}
\times\{0,\ldots,\lfloor\frac{N_2-p_2}{P_2}\rfloor-1\}$,
\begin{equation}
\Pi_{n_1,n_2}^{p_1,p_2}=(\pi_{P_1 n_1+p_1+p_1',P_2n_2+p_2+p_2'})_{0\le p_1'<P_1,0\le p_2'<P_2}
\end{equation}
and
\begin{align}
\mathbb{B} = \{(n_1,n_2)\in \NN^2\mid\; &
0\le n_1 <p_1\;\text{or}\;
0 \le n_2 < p_2\;\nonumber\\
& \text{or}\;
P_1\lfloor\frac{N_1-p_1}{P_1}\rfloor \le n_1 < N_1\;\nonumber\\
&\text{or}\;
P_2\lfloor\frac{N_2-p_2}{P_2}\rfloor \le n_2 < N_2\}.
\end{align}
It is then clear that \eqref{eq:pibord} holds since
the variables $\pi_{n_1,n_2}$ with $(n_1,n_2)\in \mathbb{B}$
are not elements of the matrices $\Pi_{n_1,n_2}^{p_1,p_2}$
with $n_1 \in \{0,\ldots,\lfloor\frac{N_1-p_1}{P_1}\rfloor-1\}$
and $n_2 \in \{0,\ldots,\lfloor\frac{N_2-p_2}{P_2}\rfloor-1\}$.
In addition, since it has been assumed that $\trace(\hh^\top \vv) = 0$ and $\Vert\hh\Vert_{\rm F} = \Vert\vv\Vert_{\rm F} = 1$,
the matrices $\Pi_{n_1,n_2}^{p_1,p_2}$ and $Y_{n_1,n_2}^{p_1,p_2}$
can be decomposed in an orthogonal manner as follows:
\begin{equation}
\begin{cases}
 \Pi_{n_1,n_2}^{p_1,p_2} = \beta_{n_1,n_2}^{p_1,p_2} \hh +  \kappa_{n_1,n_2}^{p_1,p_2}\vv + (\Pi_{n_1,n_2}^{p_1,p_2})^{\perp}\\
 Y_{n_1,n_2}^{p_1,p_2} = h_{n_1,n_2}^{p_1,p_2} \hh +  v_{n_1,n_2}^{p_1,p_2}\vv + (Y_{n_1,n_2}^{p_1,p_2})^{\perp}\\
\end{cases}
\end{equation}
where
\begin{align}
\beta_{n_1,n_2}^{p_1,p_2} &= \trace(\hh^\top \Pi_{n_1,n_2}^{p_1,p_2}),\nonumber \\
\kappa_{n_1,n_2}^{p_1,p_2} &= \trace(\vv^\top \Pi_{n_1,n_2}^{p_1,p_2}),\\
(\Pi_{n_1,n_2}^{p_1,p_2})^{\perp} &= \Pi_{n_1,n_2}^{p_1,p_2} - \beta_{n_1,n_2}^{p_1,p_2} \hh - \kappa_{n_1,n_2}^{p_1,p_2} \vv,\nonumber \\
 (Y_{n_1,n_2}^{p_1,p_2})^{\perp} &= Y_{n_1,n_2}^{p_1,p_2} - h_{n_1,n_2}^{p_1,p_2} \hh - v_{n_1,n_2}^{p_1,p_2} \vv,
\end{align}
and $(h_{n_1,n_2}^{p_1,p_2},v_{n_1,n_2}^{p_1,p_2})$ is given by \eqref{eq:yv}.
After some simplications, we have thus to minimize,
for every $n_1 \in \{0,\ldots,\lfloor\frac{N_1-p_1}{P_1}\rfloor-1\}$
and $n_2 \in \{0,\ldots,\lfloor\frac{N_2-p_2}{P_2}\rfloor-1\}$,
\begin{align}
\frac{1}{2}&\Vert \Pi_{n_1,n_2}^{p_1,p_2} - Y_{n_1,n_2}^{p_1,p_2}\Vert^2_{\rm F} + \mu\; \rho_{\rm tv}\big( \trace(\hh^\top \Pi_{n_1,n_2}^{p_1,p_2}),\trace(\vv^\top \Pi_{n_1,n_2}^{p_1,p_2})\big)\nonumber\\
&= \frac{1}{2}\Vert (\Pi_{n_1,n_2}^{p_1,p_2})^{\perp} - (Y_{n_1,n_2}^{p_1,p_2})^{\perp}\Vert^2_{\rm F} \nonumber\\
&+ \frac{1}{2}(\beta_{n_1,n_2}^{p_1,p_2} - h_{n_1,n_2}^{p_1,p_2})^2\nonumber\\
& + \frac{1}{2}(\kappa_{n_1,n_2}^{p_1,p_2} - v_{n_1,n_2}^{p_1,p_2})^2\nonumber\\
& +  \mu \; \rho_{\rm tv}\big(\beta_{n_1,n_2}^{p_1,p_2} ,\kappa_{n_1,n_2}^{p_1,p_2}\big).
\end{align}
This shows that \eqref{eq:proxtv2} is satisfied and that $(\Pi_{n_1,n_2}^{p_1,p_2})^{\perp} =  (Y_{n_1,n_2}^{p_1,p_2})^{\perp}.$
Eq. \eqref{eq:proxtv1} straightforwardly follows.

\section{Derivation of Algorithm \ref{algo:accPPA}} \label{algop:PPA_acc}
Let $\Pi_F$ denote the projector on $\ran F$. For every $u\in \RR^K$, we have
\begin{equation}
u = \Pi_F u + u^{\perp}
\label{eq:acc1}
\end{equation}
where $u^{\perp}\in (\ran F)^{\perp}$ is the projection error and 
there exists $q\in \RR^N$ such that
\begin{equation}
\Pi_F u = F q.
\label{eq:acc2}
\end{equation}
By combining this with the fact that $F^{\top}u^\perp = 0$, we obtain the relation,
\begin{equation}
q = \frac{1}{\nu} F^{\top} 
\label{eq:acc4}
\end{equation}
which allows us to deduce from \eqref{eq:acc1} that \begin{equation*}
u^{\perp} = u -\frac{1}{\nu} F F^{\top} u.                                                   
                                                    \end{equation*}
Now, consider the first step of Algorithm \ref{algo:PPA_frame}: $(\forall j\in \{1,\ldots,\JJ\})$
\begin{equation}
p_{j,\ell} = u_{j,\ell} + \frac{F}{\nu} (\prox_{\nu\gamma g_j/\omega_j} ( F^{\top}u_{j,\ell}) - F^{\top} u_{j,\ell})+a_{j,\ell}
\label{eq:acc6}
\end{equation}
where $a_{j,\ell}$ is assumed to belong to $\ran F$, i.e. $a_{j,\ell}
= F \widetilde{a}_{j,\ell}$ with $\widetilde{a}_{j,\ell}\in \RR^N$.
Defining $q_{j,\ell}\in \RR^N$ similarly to \eqref{eq:acc2} yields $\Pi_F p_{j,\ell} = F q_{j,\ell} $.
According to \eqref{eq:acc4}, $q_{j,\ell}$ is such that
\begin{equation}
q_{j,\ell} = \frac{1}{\nu} F^{\top} p_{j,\ell}.
\label{eq:acc7}
\end{equation}

\noindent By combining \eqref{eq:acc6} and \eqref{eq:acc7}, 
\begin{equation}
q_{j,\ell} = \frac{1}{\nu} \prox_{\nu\gamma g_j/\omega_j} ( v_{j,\ell}) 
+\widetilde{a}_{j,\ell}\qquad \mbox{where $v_{j,\ell} = F^{\top}u_{j,\ell}$}.
\label{eq:acc7p}
\end{equation}
Moreover, since $p_{j,\ell} = F q_{j,\ell} + p_{j,\ell}^{\perp}$, the  computation of the variable $p_\ell = \sum_{j=1}^{J} \omega_j p_{j,\ell}$ in 
Algorithm \ref{algo:PPA_frame} can be rewritten as

\begin{equation}
p_\ell = F\sum_{j=1}^{\JJ} \omega_j q_{j,\ell} + \sum_{j=1}^{\JJ} \omega_j p_{j,\ell}^{\perp} + \sum_{j=\JJ + 1}^{J} \omega_j p_{j,\ell}
\end{equation}
where, according to \eqref{eq:acc1}, \eqref{eq:acc2}, \eqref{eq:acc6}
and \eqref{eq:acc7p},
 
\begin{equation}
p_{j,\ell}^{\perp} = u_{j,\ell} - \frac{1}{\nu} F F^{\top} u_{j,\ell} = u_{j,\ell}^{\perp}.
\label{eq:acc8}
\end{equation}
In the new formulation, the last steps of the algorithm consist of updating 
$u_{j,\ell}^{\perp}$ and $v_{j,\ell}$, for all
$j \in \{1,\ldots, \JJ\}$.
We propose to define
$r_\ell = 2p_\ell - x_\ell$ , $\widetilde{r}_\ell = F^{\top} r_\ell$ and $r_\ell^{\perp} = r_\ell -\frac{1}{\nu} F \tilde{r}_\ell$,
 which yields $v_{j,\ell+1} = v_{j,\ell} + \lambda_\ell \big(\widetilde{r}_\ell - F^{\top} p_{j,\ell} \big)$ and $u_{j,\ell+1}^{\perp} = u_{j,\ell}^{\perp} + \lambda_\ell \big(r_\ell^{\perp} - p_{j,\ell}^{\perp} \big)$.
By using \eqref{eq:acc7} and \eqref{eq:acc8}, these relations 
can be simplified as

\begin{equation}
\begin{cases}
v_{j,\ell+1} = v_{j,\ell} + \lambda_\ell \big(\widetilde{r}_\ell - \nu q_{j,\ell} \big) \\
\mbox{and}\\
u_{j,\ell+1}^{\perp} = u_{j,\ell}^{\perp} + \lambda_\ell \big(r_\ell^{\perp} - u_{j,\ell}^{\perp} \big),
\end{cases}
\end{equation}
which leads to Algorithm \ref{algo:accPPA}.\\
We finally note that Assumption \ref{as:errorterms} \ref{a:3}) implies
that Assumption \ref{ass:conv_spingarn}~\ref{a:4}) is satisfied since, for every $j\in\{1,\ldots,\JJ\}$,
\begin{equation}
\sum_{\ell\in\NN}\lambda_\ell \;\Vert a_{j,\ell}\Vert
= \sum_{\ell\in\NN}\lambda_\ell \;\Vert F \widetilde{a}_{j,\ell}\Vert 
\le \|F\| \sum_{\ell\in\NN}\lambda_\ell \;\Vert \widetilde{a}_{j,\ell}\Vert<+\infty.
\end{equation}
This allows us to transpose the convergence results concerning 
Algorithm \ref{algo:PPA_frame} to Algorithm  \ref{algo:accPPA}.

\end{document}